\documentclass[12pt,titlepage,a4paper]{amsart}
\usepackage{graphicx}
\usepackage{amssymb}
\usepackage{amsfonts}
\usepackage{epsfig}

\def\SetTableau#1#2#3#4{%
  \gdef\Tabvrule{\vrule\vrule width-0.4pt}
  \gdef\Tabhrule{\hrule\hrule height-0.4pt}  
  \gdef\Tabstrut{\vrule height#1 depth#2 width0pt\relax}
  \gdef\Tabbox##1{\hbox to #3{\hskip0.4pt\hfill\Tabstrut$#4##1$\hfill}}
} 
\def\PetitTableau{\SetTableau{1.65ex}{0.55ex}{2.2ex}{\scriptstyle}}

\def\Case#1{\vcenter{\Tabhrule%
                   \hbox{\Tabvrule\Tabbox{#1}\Tabvrule}\Tabhrule}}

\def\GenTab#1{\vcenter{\halign{&$\Case{##}$\cr#1}}\egroup}
\def\Tableau{%
  \bgroup%
  \let\ =\omit%
  \let\\=\cr%
  \offinterlineskip\GenTab}
\newtheorem{example}{Example}[section]

\newtheorem{theorem}[example]{Theorem}

\newtheorem{corollary}[example]{Corollary}
\newtheorem{conjecture}[example]{Conjecture}
\newtheorem{definition}[example]{Definition}
\newtheorem{proposition}[example]{Proposition}
\newtheorem{algorithm}[example]{Algorithm}
\newtheorem{lemma}[example]{Lemma}

\def\Proof{\noindent \it Proof -- \rm}
\def\qed{\hspace{3.5mm} \hfill \vbox{\hrule height 3pt depth 2 pt width 2mm}
\bigskip}

\def\QSym{{\it QSym}}          
\def\NCSF{{\bf Sym}}           
\def\FQSym{{\bf FQSym}}        
\def\MQSym{{\bf MQSym}}        
\def\FSym{{\bf FSym}}          
\def\PBT{{\bf PBT}}            
\def\PQSym{{\bf PQSym}}        
\def\CQSym{{\bf CQSym}}        
\def\SQSym{{\bf SQSym}}        
\def\SCQSym{{\bf SCQSym}}      

\def\J{{\bf J}}
\def\P{{\bf P}}                 %
\def\C{{\mathbb C}}

\def\a{{\bf a}}
\def\b{{\bf b}}
\def\c{{\bf c}}
\def\d{{\bf d}}

\def\pack{{\rm pack}}          
\def\Lie{{\rm Lie}}            
\def\L{{\mathfrak L}}
\def\ncbinomial#1#2{\left[\,\begin{matrix}#1 \cr #2\end{matrix}\,\right]}
\def\binomial#1#2{\left(\,\begin{matrix}#1 \cr #2\end{matrix}\,\right)}

\def\last{{\rm last}}   
\def\ev{{\rm Ev}}       
\def\congru{\equiv}     
\def\ssh{\Cup}          
\def\saug{\uplus}       
\def\sconc{\bullet}     
\def\Std{{\rm Std}}     
\def\Park{{\rm Park}}   
\def\convol{{*}}        
\def\convP{{*_P}}       
\def\<{\langle}
\def\>{\rangle}
\def\NN{{\mathbb N}}    
\def\RR{{\mathbb R}}    
\def\CC{{\mathbb C}}    

\def\park{{\bf a}} 

\def\F{{\bf F}}         
\def\S{{\bf S}}         
\def\G{{\bf G}}         
\def\SG{{\mathfrak S}}  

\def\dim{{\rm dim}}

\def\ch{\operatorname {ch}}

\def\End{\operatorname{End}} 

\def\PF{{\rm PF}}   
\def\PPF{{\rm PPF}} 
\def\PW{{\rm PW}}   

\def\shuff#1#2{\mathbin{
\hbox{\vbox{ \hbox{\vrule \hskip#2 \vrule height#1 width 0pt
}%
\hrule}%
\vbox{ \hbox{\vrule \hskip#2 \vrule height#1 width 0pt
\vrule }%
\hrule}%
}}}

\def\shuf{{\mathchoice{\shuff{7pt}{3.5pt}}%
{\shuff{6pt}{3pt}}%
{\shuff{4pt}{2pt}}%
{\shuff{3pt}{1.5pt}}}}%
\def\shuffle{\,\shuf\,}

\def\detassmax{unp}   
\def\PQR{{\rm PQR}}   
\def\PQS{{\rm PQS}}   

\def\raff{\succeq}    
\def\assh{\uplus}     

\def\gaudend{\ll}     
\def\droitdend{\gg}   
\def\gautrid{\!\prec\!}   
\def\miltrid{\circ}       
\def\droittrid{\!\succ\!} 

\def\qrpark{{\bf q}}  
\def\qrparkb{{\bf r}} 
\def\LPQ{{\bf LPQ}}   
\def\NC{{\rm NC}}     
\def\park{{\bf a}}    
\def\parkb{{\bf b}}   
\def\vp{{\rm vp}}     

\def\S{{\bf S}}       
\def\T{{\bf T}}       

\def\M{{\bf M}}       
\def\V{{\bf V}}       
\def\MA{{\mathcal M}} 

\def\PCat{{\bf P}}    
\def\RCat{{\bf R}}    
\def\MM{{\mathcal M}} 
\def\TD{{\mathfrak{TD}}}
\def\PS{{\bf P}}      
\def\RS{{\bf R}}      
\def\QS{{\bf Q}}      
\def\PSC{{\it P}}     
\def\RSC{{\it R}}     
\def\QSC{{\it Q}}     
\def\FSC{{\it F}}     

\def\TT{{\mathcal T}} 
\def\tasse{{\it t}}   
\def\PP{{\mathcal P}}
\def\Ig{{\bf I}}      

\def\sep{\,|\,}       
\def\sepb{|}          
\def\psev{{\rm ps}}   
\def\limproj{{\rm proj\,lim}}
\def\K{{\mathbb K}}   
\def\Sym{{\bf Sym}}   

\title[Hopf algebras and parking functions]{Hopf algebras and dendriform
structures \\ arising from parking functions}

\author[J.-C.~Novelli and J.-Y.~Thibon]
{Jean-Christophe Novelli and Jean-Yves Thibon}

\address[] {Institut Gaspard Monge, Universit\'e de Marne-la-Vall\'ee \\
5 Boulevard Descartes \\Champs-sur-Marne \\77454 Marne-la-Vall\'ee cedex 2 \\
FRANCE}
\email[Jean-Christophe Novelli]{novelli@univ-mlv.fr}
\email[Jean-Yves Thibon]{jyt@univ-mlv.fr} 
\date{}

\begin{document}

\begin{abstract}
We introduce a graded  Hopf algebra based on the set of parking
functions (hence of dimension $(n+1)^{n-1}$ in degree $n$). 
This algebra can be embedded into a noncommutative polynomial algebra
in infinitely many variables.
We determine its structure, and show that it admits natural quotients
and subalgebras whose graded components have dimensions respectively given
by the Schr\"oder numbers (plane trees), the Catalan numbers, and powers of 3.
These smaller algebras are always bialgebras and belong to some family of di-
or tri-algebras occuring in the works of Loday and Ronco.

Moreover, the fundamental notion of parkization allows one to endow
the set of parking functions of fixed length with an associative
multiplication (different from the one coming from the Shi arrangement),
leading to a generalization of the internal product of symmetric functions.
Several of the intermediate algebras are stable under this operation.
Among them, one  finds the Solomon descent algebra but also a new algebra
based on a Catalan set, admitting the Solomon algebra as a left ideal.
\end{abstract}

\maketitle
{\footnotesize
\tableofcontents
}

\newpage
\section{Introduction}

Many examples of graded Hopf algebras based on combinatorial structures occur
in apparently remote contexts. One of them is the theory of operads.  It is
quite common there that in a given operad, the free algebra on one generator
admit a Hopf structure \cite{Holt}. This structure often has an elegant
combinatorial description, the best known example being the free dendriform
algebra on one generator, also known as the Loday-Ronco algebra of planar
binary trees \cite{LR1,HNT}.

On another hand, such Hopf algebras also occur in the theory of noncommutative
symmetric functions \cite{NCSF1}, for which one central problem is to
understand complicated commutative formulas by means of simpler non
commutative analogues. It has been found over the years that such an
understanding required the introduction of larger and larger Hopf algebras,
based on more and more complex combinatorial objects.  For such algebras to be
useful in this context, it is necessary that their elements can be realized as
polynomials in some auxiliary infinite set of variables (commutative or not),
so as to recover ordinary symmetric functions after a chain of standard
manipulations (such as imposing commutation relations among the variables or
taking sums to reestablish complete symmetry). The best illustration of this
approach is provided by the algebra of Free Quasi-Symmetric Functions $\FQSym$
\cite{NCSF6}. This is an algebra of noncommutative polynomials $\F_\sigma(A)$
labeled by permutations. It contains a subalgebra $\FSym$ spanned by free
Schur functions $\S_t(A)$, labeled by standard Young tableaux. This
observation essentially amounts to a one-line proof of the
Littlewood-Richardson rule. Abstractly, however, $\FQSym$ and $\FSym$ are
isomophic to the Hopf algebras previously introduced in \cite{MR} and
\cite{PR}, and it is the polynomial realization which allows such a direct
application to symmetric functions.

Interestingly, it is the very same realization which allowed a new
understanding of the algebra $\PBT$ of planar binary trees \cite{HNT}. It
could be put on the same footing as $\FSym$, using the sylvester
correspondence instead of Robinson-Schensted, so that both algebras appear now
as special cases of a general construction.

The aim of the present article is to introduce a new extension of $\FQSym$,
that is, a larger Hopf algebra built from the same principles, but leaving
enough room to accomodate several new combinatorial Hopf algebras.

It turns out that most of the Hopf algebras arising in the process also have
an operadic interpretation, in general as some kind of trialgebra or dialgebra
\cite{Lod-dend,LRtri}, thus providing polynomial realizations of those as
well.

Our master algebra, denoted by $\PQSym$, for \emph{Parking Quasi-Symmetric
Functions}, is built on the set of parking functions, a special family of
words which can in many respects be regarded as natural generalizations of
permutations. Geometrically, permutations correpond to chambers of the
Coxeter arrangement of type $A_{n-1}$, while parking functions label those of
the Shi arrangement \cite{Ath}, but this is not the only possible explanation
(see, \emph{e.g.}, \cite{LodPark}), and our choice was rather dictated by
elementary combinatorial considerations (see Appendix).

Our first task will be to elucidate the structure of $\PQSym$.
It will be shown that it is free, cofree, and actually self-dual, with a free
primitive Lie algebra.
This will be done by means of Foissy's theory of bidendriform
bialgebras \cite{Foi}. Next, we shall determine explicit generators
and multiplicative bases of $\PQSym$ and $\PQSym^*$. Then come
the realizations, given by simple and explicit noncommutative polynomials
for the natural basis of $\PQSym^*$, and in terms of integer matrices,
reminescent of the construction of $\MQSym$ \cite{NCSF6}, for the natural
basis of $\PQSym$ itself. After that, we shall start the investigation of
smaller Hopf algebras arising from $\PQSym$ by natural processes.

Recall that the dimension of $\PQSym$ in degree $n$ is $(n+1)^{n-1}$. We shall
show that it admits natural quotients and subalgebras whose graded components
have dimensions respectively given by the Schr\"oder numbers (plane trees),
the Catalan numbers, powers of 3 and powers of 2.
Most of those turn out to be related to the theory of operads, and to belong
to some family of di- or tri-algebras occuring in the works of Loday and
Ronco. 
We shall in particular recover the free dendriform trialgebra on one generator
(Schr\"oder numbers) and the free cubical trialgebra.
Similarly, we obtain a cocommutative Hopf algebra based on a Catalan
set, which is isomorphic to the free dendriform dialgebra on one generator
as an algebra, but not as a coalgebra.

Moreover, the fundamental notion of parkization of a word, which is needed
from the beginning,  allows one to endow the set of parking functions of fixed
length with an associative multiplication (different from the one coming from
their interpretation as chambers of the Shi arrangement), leading to a
generalization of the internal product of symmetric functions.
Several of the intermediate algebras are stable under this operation.
Among them, one finds the Solomon descent algebra and the Solomon-Tits
algebra, but also a new algebra based on a Catalan set, admitting the Solomon
algebra as a left ideal.

\smallskip
This paper is structured as follows: the preliminaries present some background
about parking functions and dendriform structures needed in the sequel and
give a realization the free dendriform trialgebra on one generator in terms of
noncommutative polynomials.
In Section~\ref{hopfalg}, we present our principal algebra $\PQSym$, and
investigate its most important features, mostly reying upon its bidendriform
bialgebra structure. We then move to a subalgebra $\SQSym$ of $\PQSym$, whose
Hilbert series is given by the little Schr\"oder numbers and prove in
particular that it is isomorphic to the free dendriform trialgebra on one
generator (Section~\ref{sqsym}).
In Section~\ref{cqsym}, we study another subalgebra $\CQSym$ of $\PQSym$ whose
Hilbert series is given by the Catalan numbers, show that it is cocommutative,
that it is stable under the internal product of $\PQSym$ and that its dual is
a natural generalization of $\QSym$.
In Section~\ref{cc}, we present $\SCQSym$, a quotient of $\SQSym$ whose
Hilbert series is given by powers of $3$ and show in particular that it is
isomorphic to the free cubical trialgebra on one generator.
Finally, the Appendix presents how the construction of $\PQSym$ arose from
considerations about free probability and an exercise proposed by Kerov in
1995.
Most of these results were announced in~\cite{NT1}.

\medskip
{\footnotesize
{\it Acknowledgements.-}
This project has been partially supported by CNRS and by EC's IHRP Programme,
grant HPRN-CT-2001-00272, ``Algebraic Combinatorics in Europe".
The authors would also like to thank the contributors of the MuPAD project,
and especially those of the combinat package, for providing the development
environment for this research (see~\cite{HT} for an introduction to
MuPAD-Combinat).
}

\newpage
\section{Preliminaries}

\subsection{Notations}

Our notations for ordinary symmetric functions will be those of \cite{Mcd}.
Other undefined  notations can be found in \cite{NCSF1,NCSF6}, although the
essential ones will be recalled when needed.

\subsubsection{}
To start with, we shall need  the following two operations on words.

For a word $w$ on the alphabet $\{1,2,\ldots\}$, denote by $w[k]$ the word
obtained by replacing each letter $i$ by the integer $i+k$.
If $u$ and $v$ are two words, with $u$ of length $k$, one defines
the {\em shifted concatenation}
\begin{equation}
u\sconc v = u\cdot (v[k])
\end{equation}
and the {\em shifted shuffle}
\begin{equation}
u\ssh v= u\shuffle (v[k])\,.
\end{equation}
where $\shuffle$ is the usual shuffle product on words defined by
\begin{equation}
(au)\shuffle (bv) = a\cdot(u\shuffle (bv)) + b\cdot ((au)\shuffle v),
\end{equation}
with $u\shuffle\epsilon=\epsilon\shuffle u=u$ if $\epsilon$ is the empty word.

It is immediate to see that the set of permutations is closed under both
operations. The subalgebra spanned by those elements is isomorphic to
the convolution algebra of symmetric groups \cite{MR} or to Free
Quasi-Symmetric Functions \cite{NCSF6}, whose definition is recalled below.

\subsubsection{}
Let $A$ be a totally ordered alphabet. We denote by $\K$ a field
of characteristic $0$, and by $\K\<A\>$ the free associative algebra
over $A$ when $A$ is finite, and the projective limit
$\limproj_B \K\<B\>$, where $B$ runs over finite subsets of $A$,
when $A$ is infinite, which will be generally assumed in the sequel.

Given a totally ordered alphabet $A$, the \emph{evaluation vector} $\ev(w)$ of
a word $w$ is the sequence of the numbers of occurrences of all the elements
of $A$.

Recall that the standardized $\Std(w)$ of a word $w\in A^*$ is
the permutation obtained by iteratively scanning $w$ from
left to right, and labelling $1,2,\ldots$ the occurrences of its
smallest letter, then numbering the occurrences of the next one, and
so on. Alternatively, $\sigma=\Std(w)^{-1}$ can be characterized as
the unique permutation of minimal length such that $w\sigma$ is a
nondecreasing word. For example, $\Std(bbacab)=341624$.

This characterizes completely the sequences of transpositions
effected by the bubble sort algorithm on $w$. An elementary observation,
which is at the basis of the constructions of \cite{NCSF6}, is that
the noncommutative polynomials
\begin{equation}
\G_\sigma(A) =\sum_{w\in A^*;\Std(w)=\sigma}w
\end{equation}
span a subalgebra of $\K\langle A\rangle$. Moreover, if $A$ is
infinite, this subalgebra admits a natural Hopf algebra structure.
This is $\FQSym$, the algebra of {\em Free Quasi-Symmetric Functions}.

Let $\F_\sigma=\G_{\sigma^{-1}}$. The coproduct is defined by
\begin{equation}
\Delta\F_\sigma=\sum_{u\cdot v=\sigma}\F_{\Std(u)}\otimes \F_{\Std(v)}\,,
\end{equation}
where $u\cdot v$ means concatenation. The scalar product is defined by
\begin{equation}
\<\F_\sigma\,,\,\G_\tau\>=\delta_{\sigma,\tau}\,,
\end{equation}
where $\delta$ is the Kronecker symbol,
and one has then for all $F,G,H\in\FQSym$
\begin{equation}
\<FG,H \>=\<F\otimes G,\Delta H\>\,.
\end{equation}
The product formula in the $\F$ basis is
\begin{equation}
\F_\alpha\F_\beta=\sum_{\gamma\in\alpha\ssh\beta }\F_\gamma\,.
\end{equation}
The sum of the inverses of the permutations occuring in
$\alpha^{-1}\ssh\beta^{-1}$ is called \emph{convolution} and denoted by
$\alpha * \beta$ \cite{Re,MR}.

\subsubsection{}
A general process for constructing interesting subalgebras of $\FQSym$ is to
take sums of the form
\begin{equation}
\PP_{\bf x}(A) = \sum_{\mathfrak{P}(\sigma)={\bf x}} \F_\sigma\,,
\end{equation}
where $\mathfrak{P}$ is the left symbol of some Robinson-Schensted type
correspondence. If we take the original Robinson-Schensted map, we obtain
$\FSym$, the algebra of free symmetric functions~cite{NCSF6}. If we take the
sylvester congruence~\cite{HNT}, we obtain $\PBT$, the Loday-Ronco algebra of
planar binary trees. Finally, if we take the hypoplactic correspondence
\cite{NCSF4}, we obtain $\Sym$, the algebra of noncommutative symmetric
functions. The dual Hopf algebras are obtained in each case by imposing the
corresponding congruence (plactic, sylvester, hypoplactic) on $A^*$.

\subsection{Parking functions}

In the following, we shall see that it is possible to replace
permutations by {\em parking functions} is all these constructions.
As one will see, it is obvious that the set of parking functions
is stable under shifted concatenation and shifted shuffle, and
many other classes of words share this property. The point is that
for parking functions, the resulting algebra has a natural
Hopf structure, and that it is again possible to find a polynomial
realization. Moreover, an interesting internal product can be defined.

\subsubsection{}

A \emph{parking function} is a word $\park=a_1a_2\cdots a_n$ of length $n$ on
$[n]=\{1,2,\ldots,n\}$ whose nondecreasing rearrangement
$\park^\uparrow=a'_1a'_2\cdots a'_n$ satisfies $a'_i\le i$ for all $i$.
Let $\PF_n$ be the set of such words.

For example, $\PF_1= \{1\}$, $\PF_2 = \{11, 12, 21\}$, and

\begin{equation}
\begin{split}
\PF_3 = \{&111, 112, 121, 211, 113, 131, 311, 122, 212, 221,\\
& 123, 132, 213,
231, 312, 321\}
\end{split}
\end{equation}

\subsubsection{}
It is well-known that $|\PF_n|=(n+1)^{n-1}$, and that the
permutation representation of $\SG_n$ naturally supported  by $\PF_n$
has Frobenius characteristic (see~\cite{Haiman})
\begin{equation}\label{valgn}   
(-1)^n\omega(h_n^*) 
\end{equation}
where $f\mapsto f^*$ is the involution on symmetric functions defined 
on the genertors $h_n$ as follows
(see~\cite{Mcd}, ex.~24 p.~35).
If we set $H(t):=\sum_{n\ge 0}h_nt^n$ and $H^*(u):=\sum_{n\ge 0}h_n^*u^n$,
then
\begin{equation}
u=tH(t) \ \Longleftrightarrow \ t=uH^*(u).
\end{equation}

Each nondecreasing parking function generates a sub-permutation representation
of $\PF_n$.
It is easy to see that the number of nondecreasing parking functions of length
$n$ is the Catalan number $C_n=\frac{1}{n+1}\binom{2n}{n}$.

\subsubsection{Prime parking functions}

This important notion has been derived by Gessel in 1997 (see~\cite{Stan2}). 
Given a parking function of length $n$, one says that
$b\in\{0,1,\ldots,n\}$ is a \emph{breakpoint} of $\park$ if
$|\{i\,|\, a_i\le b\}|=b$.
For example, the parking function $112256679$ has the five breakpoints
$\{0,4,5,8,9\}$.
Then, $\park\in \PF_n$ is said to be \emph{prime} if its only breakpoints are
the trivial ones: $0$ and $n$. 
Let $\PPF_n\subset\PF_n$ be the set of prime parking functions on $[n]$.
For example, 
\begin{equation}
\PPF_1=\{1\}, \quad
\PPF_2=\{11\}, \quad
\PPF_3=\{111,112, 121, 211\}.
\end{equation}
It can easily be shown that $|\PPF_n|=(n-1)^{n-1}$ for $n\geq2$
(see~\cite{Stan2,Kal} and Section~\ref{compos-sec}).
The number of nondecreasing prime parking functions of length $n$
is the shifted Catalan number $C_{n-1}$: they are obtained by concatenating a $1$ to the
left of all nondecreasing parking functions of length $n-1$.

As already mentioned, it is immediate to see that the set of all parking
functions is closed under shifted concatenation and shifted shuffle. The
prime parking functions are exactly those that do not occur in any nontrivial
shifted shuffle of parking functions.
This observation is at the basis of our definition of the Hopf algebra of
parking functions (see Section~\ref{hopfalg}).

\subsubsection{The module of prime parking functions}

Parking functions can be classified according to the factorization of their
nondecreasing reorderings $\park^\uparrow$ with respect to the operation of
shifted concatenation. That is, if
\begin{equation}
\park^\uparrow = w_1 \sconc w_2 \sconc\cdots\sconc w_r
\end{equation}
is the unique maximal factorization of $\park^\uparrow$, each $w_i$ is a
nondecreasing prime parking function.
Let us define $i_k=|w_k|$ and let $I=(i_1,\ldots,i_r)$. 
We shall say that $\park$ is of \emph{type} $I$ and denote by $\PPF_I$ the set
of parking functions of type $I$.
For example, the parking function $966142272$ is of type $(1,4,3,1)$ and
the number of parking functions of length $4$ of each type is

$$
\begin{array}{|c|c|c|c|c|c|c|c|}
\hline
(4) & (31) & (13) & (22) & (211) & (121) & (112) & (1111) \\
\hline
27  & 16   &  16  &  6   &   12  &  12   &  12   &   24   \\
\hline
\end{array}
$$

\bigskip
The set $\PPF_{n}$ of prime parking functions of length $n$ obviously
is a sub-permutation representation of $\PF_n$.
It can be shown that its Frobenius characteristic is 
\begin{equation}\label{valfn}
f_n=-\omega(e_n^*)
\end{equation}
(see the Appendix for a direct proof).
One can also obtain it as follows.

The set $\PPF_I$ of parking functions of type $I$ is a
sub-permutation representation of $\PF_n$, and its Frobenius
characteristic is
\begin{equation}
\ch(\PPF_I) = f_{i_1} \ldots f_{i_r}\,,
\end{equation}
since it is induced from the permutation representation
of the Young subsgroup $\SG_I$
on the Cartesian product $\PPF_{i_1}\times\cdots\times\PPF_{i_r}$.
Now, $\PF_n=\bigsqcup_{I\vDash n}\PPF_I$, so that
\begin{equation}
g_n = \sum_{I\vDash n}f_{i_1} \ldots f_{i_r}\,,
\end{equation}
which amonts to
\begin{equation}
g :=\sum_{n\ge 0}g_n = (1-f)^{-1}\,\ \text{where} \ f=\sum_{n\ge 1}f_n\,.
\end{equation}
Thus, if we know that $g_n$ is given by (\ref{valgn}), we obtain that $f_n$ is
given by (\ref{valfn}), and conversely.
A noncommutative version of these results will be established in
Section~\ref{compos-sec}.

\subsection{Dendriform dialgebras}

A \emph{dendriform dialgebra}, as defined by Loday~\cite{Lod-dend}, is an
associative algebra $D$ whose multiplication $\odot$ splits into two binary
operations
\begin{equation} 
x \odot y = x\gaudend y + x\droitdend y\,,
\end{equation}
called left and right, satisfying the  following three compatibility relations
for all $a$, $b$, and $c$ different from $1$ in $D$:
\begin{equation}
\label{dend1}
(a\gaudend b)\gaudend c = a\gaudend (b\odot c)
\end{equation}
\begin{equation}
\label{dend2}
(a\droitdend b)\gaudend  c = a\droitdend  (b\gaudend c)
\end{equation}
\begin{equation}
\label{dend3}
(a\odot b)\droitdend  c = a\droitdend  (b\droitdend c)
\end{equation}

These relations are satisfied by shuffle algebras with $\odot=\shuffle$ and
for $x=ua$ and $y=vb$ ($a$ and $b\in A$),
\begin{equation}
x\droitdend y = (ua\shuffle v) b\,,\quad
x\gaudend y  = (u\shuffle vb) a
\end{equation}

It turns out that the free associative algebra $\K\langle A\rangle$ is also a
dendriform dialgebra. Actually, it is even a dendriform trialgebra, as
explained below.

\subsection{Dendriform trialgebras}
\label{trigd1}

A \emph{dendriform trialgebra}~\cite{LRtri} is an associative algebra whose
multiplication $\odot$ splits into three pieces
\begin{equation}
x\odot y = x\gautrid y + x\miltrid y + x\droittrid y\,,
\end{equation}
where $\miltrid$ is associative, and
\begin{eqnarray}
(x\gautrid y)\gautrid z = x\gautrid (y\odot z)\,,\\
(x\droittrid y)\gautrid z = x\droittrid (y\gautrid z)\,,\\
(x\odot y)\droittrid z = x\droittrid (y\droittrid z)\,,\\
(x\droittrid y)\miltrid z = x\droittrid (y\miltrid z)\,,\\
(x\gautrid y)\miltrid z = x\miltrid (y\droittrid z)\,,\\
(x\miltrid y)\gautrid z = x\miltrid (y\gautrid z)\,.
\end{eqnarray}

Let $A=\{a_1<a_2<\cdots <a_n <\cdots \,\}$ be an infinite linearly ordered
alphabet. Recall that $\K\langle A\rangle$ is understood as the projective
limit of the $\K\langle A_n\rangle$ where $A_n$ is the interval $[a_1,a_n]$ of
$A$.
We denote by $\max(w)$ the greatest letter occuring in the word $w\in A^*$.

\medskip
\begin{definition}
\label{def-trigdend}
For two non empty words $u,v\in A^*$, we set
\begin{eqnarray}
u\gautrid v=\begin{cases} uv &
\text{if $\max(u)>\max(v)$}\cr 0 &\mbox{otherwise.} \end{cases}\\
u\miltrid v=\begin{cases} uv &
\text{if $\max(u)=\max(v)$}\cr 0 &\mbox{otherwise,} \end{cases}\\
u\droittrid v=\begin{cases} uv &
\text{if $\max(u)<\max(v)$}\cr 0 &\mbox{otherwise,} \end{cases}
\end{eqnarray}
\end{definition}

\begin{lemma}
\label{lem-trigdend}
The three operations $\gautrid\,,\ \miltrid\,,\ \droittrid$, endow the
augmentation ideal $\K\<A\>^+$ with the structure of a dendriform trialgebra.
\end{lemma}

\Proof A straightforward verification.
\qed

Setting $\gaudend=\gautrid$ and $\droitdend=\miltrid+\droittrid$, we obtain
a dendriform dialgebra.

It is known~\cite{LRtri} that the free dendriform trialgebra on one generator,
denoted here by $\TD$, is a free associative algebra with Hilbert series
\begin{equation}
\label{sg-dendt}
\sum_{n\geq0} s_n t^n = \frac{1+t-\sqrt{1-6t+t^2}}{4t}
= 1 + t + 3t^2 + 11t^3 + 45t^4 + 197t^5 + \cdots
\end{equation}
that is, the generating function of the \emph{super-Catalan}, or
\emph{little Schr\"oder} numbers, counting \emph{plane trees}.

The previous considerations allow us to give a simple polynomial realization
of $\TD$. Consider the polynomial
\begin{equation}
\M_1=\sum_{i\ge 1}a_i \,,
\end{equation}
(the sum of all letters). We can then state:

\begin{theorem}
The sub-trialgebra of $\K\<A\>^+$ generated by $\M_1$ is free as a dendriform
trialgebra.
\end{theorem}

We shall need the following construction on words.
With any word $w$ of length $n$, associate a plane tree $\TT(w)$ with $n+1$
leaves, as follows: if $m=\max(w)$ and if $w$ has exactly $k$ occurences of
$m$, write
\begin{equation}
w=v_0\,m\,v_1\,m\,v_2\cdots v_{k-1}\,m\,v_k\,,
\end{equation}
where the $v_i$ may be empty. Then, $\TT(w)$ is the tree obtained by
grafting the subtrees $\TT(v_0),\TT(v_1),\ldots,\TT(v_k)$ (in this order)
on a common root, with the initial condition $\TT(\epsilon)=\emptyset$
for the empty word.

\begin{figure}[ht]
\centerline{\epsfxsize=6cm \epsfbox{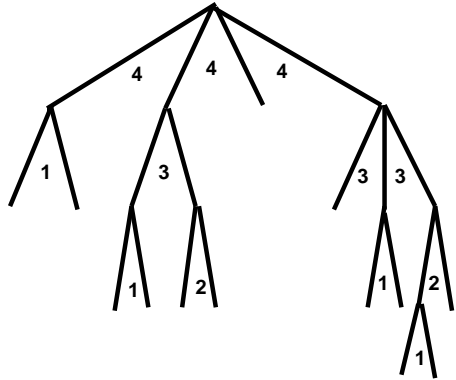}}
\caption{\label{arb-plan}The tree of $141324431312$.}
\end{figure}

For example, the tree associated with $141324431312$ is represented in
Figure~\ref{arb-plan}.

Now associate with each plane tree $T$ a polynomial by
\begin{equation}
\MM_T := \sum_{\TT(w)=T}w\,.
\end{equation}
These belong to the subtrialgebra generated by $\M_1$ since, if $T$
has as subtrees of its root $T_1,\ldots,T_k$, one has
\begin{equation}
\MM_T = \MM_{T_1}\droittrid\M_{1} \miltrid (\MM_{T_2}\droittrid\M_{1})
\miltrid \cdots \miltrid (\MM_{T_{k-1}}\droittrid\M_{1}) \gautrid \M_{T_k}.
\end{equation}
For example, with the tree $T$ presented in Figure~\ref{arb-plan}, one gets
the expression:
\begin{equation}
\M_{1}\droittrid\M_{1} \miltrid
((\M_{1}\droittrid\M_{1}\gautrid\M_{1})\droittrid\M_{1}) \miltrid
\M_{1} \gautrid \left(\M_{1}\miltrid (\M_{1}\droittrid\M_{1})\gautrid
(\M_{1}\droittrid\M_{1})\right).
\end{equation}

\Proof[of the theorem]
Since it is already known that the dimension of the free dendriform trialgebra
on one generator has dimensions given by the little Schr\"oder numbers,
we just need to show that all terms of the Hilbert series of this subalgebra
are greater than or equal to the terms of Equation~(\ref{sg-dendt}).
The polynomials $\MM_T$, being sums over disjoint sets of words, are
obviously linearly independent, whence the result.
\qed

\begin{corollary}[\cite{LRtri}]
The free commutative dendriform trialgebra on one generator
is $QSym^+$, the augmentation ideal of quasi-symmetric functions.
\end{corollary}

Indeed, it is the image of $\TD$ by the ring homomorphism mapping the letters
$a_i$ to commuting variables $x_i$.

Other applications of this realization of $\TD$ will be given in
Section~\ref{sqsym}.

\subsection{Bidendriform bialgebras}

These have been introduced by Foissy in~\cite{Foi}.
A \emph{bidendriform bialgebra} is a  dendriform dialgebra equipped with a
coproduct that splits into two parts, satifying the
\emph{codendriform relations}, obtained by dualizing the dendriform
relations, and certain compatibility properties with the two half-products.

A \emph{codendriform coalgebra} is a coalgebra $C$ whose coproduct $\Delta$
splits as $\Delta(c)=\overline\Delta(c)+c\otimes1+1\otimes c$ and
$\overline\Delta=\Delta_\gaudend  + \Delta_\droitdend $, such that, for all
$c$ in $C$:
\begin{equation}
\label{codend1}
(\Delta_\gaudend \otimes Id) \circ \Delta_\gaudend (a) =
(Id\otimes\overline\Delta)\circ\Delta_\gaudend (a),
\end{equation}
\begin{equation}
\label{codend2}
(\Delta_\droitdend \otimes Id) \circ \Delta_\gaudend (a) =
(Id\otimes\Delta_\gaudend )\circ\Delta_\droitdend (a),
\end{equation}
\begin{equation}
\label{codend3}
(\overline\Delta\otimes Id)\circ\Delta_\droitdend (a) =
(Id\otimes\Delta_\droitdend ) \circ \Delta_\droitdend (a).
\end{equation}

The Loday-Ronco algebra of planar binary trees introduced in~\cite{LR1} arises
as the free dendriform dialgebra on one generator. This is moreover a Hopf
algebra, which turns out to be self-dual, so that it is also codendriform.

There is some compatibility between the dendriform and the codendriform
structures, leading to what has been called by Foissy~\cite{Foi} a
\emph{bidendriform bialgebra}.
A bidendriform bialgebra is both a dendriform dialgebra and a codendriform
coalgebra satisfying the following four compatibility relations
\begin{equation}
\label{bidend1}
\Delta_\droitdend  (a\droitdend b) =
      a'b'_\droitdend \!\otimes\!  a''\!\!\droitdend \!b''_\droitdend \,
\,+\, a'\!\otimes\! a''\!\!\droitdend \!b \,
\,+\, b'_\droitdend \!\otimes\! a\!\droitdend \!b''_\droitdend  \,
\,+\, ab'_\droitdend \!\otimes\! b''_\droitdend  \,
\,+\, a\!\otimes\! b\,,
\end{equation}
\begin{equation}
\label{bidend2}
\Delta_\droitdend  (a\gaudend b) =
      a'b'_\droitdend \!\otimes\! a''\!\gaudend \!b''_\droitdend 
\,+\, a'\!\otimes\! a''\!\gaudend \!b 
\,+\, b'_\droitdend \!\otimes\! a\!\gaudend \!b''_\droitdend \,,
\end{equation}
\begin{equation}
\label{bidend3}
\Delta_\gaudend  (a\droitdend b) =
      a'b'_\gaudend \!\otimes a''\!\droitdend \!b''_\gaudend
\,+\, ab'_\gaudend \otimes b''_\gaudend  
\,+\, b'_\gaudend \otimes a\droitdend b''_\gaudend \,,
\end{equation}
\begin{equation}
\label{bidend4}
\Delta_\gaudend  (a\gaudend b) =  a'b'_\gaudend \!\otimes\! a''\!\gaudend
\!b''_\gaudend 
\,+\, a'b\!\otimes\! a'' 
\,+\, b'_\gaudend \!\otimes\! a\!\gaudend \!b''_\gaudend  \,+\, b\!\otimes\! a\,,
\end{equation}
where the pairs $(x',x'')$ (resp. $(x'_\gaudend ,x''_\gaudend )$ and
$(x'_\droitdend,x''_\droitdend)$) correspond to all possible elements
occuring in $\overline\Delta x$ (resp. $\Delta_\gaudend x$ and
$\Delta_\droitdend  x$), summation signs being understood (Sweedler's
notation).

Foissy has shown~\cite{Foi} that a connected bidendriform bialgebra
${\mathcal B}$ is always free as an associative algebra and self-dual as a
Hopf algebra. Moreover, its primitive Lie algebra is free, and as a dendriform
dialgebra, ${\mathcal B}$ is also free over the space of totally primitive
elements (those annihilated by $\Delta_\gaudend$ and $\Delta_\droitdend$).

It is also proved in \cite{Foi} that $\FQSym$ is bidendriform, so that it
satisfies all these properties. 

\newpage
\section{The Hopf algebra of parking functions}
\label{hopfalg}

\subsection{The algebra $\PQSym$}

Since permutations are special parking functions and  parking functions
are stable under the shifted shuffle, it is natural to embed the algebra of
Free Quasi-Symmetric functions $\FQSym$ of~\cite{NCSF6} into an
algebra spanned by  elements $\F_\park$ ($\park\in\PF$), with the same
multiplication rule:
\begin{equation}
\label{prodF}
\F_{\park'}\F_{\park''}:=\sum_{\park\in\park'\ssh\park''}\F_\park\,.
\end{equation}
We shall call this algebra $\PQSym$ (Parking Quasi-Symmetric functions).

For example,
\begin{equation}
\F_{1}\F_{1}= \F_{12} + \F_{21}, \qquad
\F_{1}\F_{11}= \F_{122} + \F_{212} + \F_{221}\,.
\end{equation}
\begin{equation}
\F_{1}\F_{12}= \F_{123} + \F_{213} + \F_{231}, \quad
\F_{1}\F_{21}= \F_{132} + \F_{312} + \F_{321}\,.
\end{equation}
\begin{equation}
\F_{12}\F_{11}= \F_{1233} + \F_{1323} + \F_{1332} + \F_{3123} + \F_{3132}
+ \F_{3312}\,.
\end{equation}
\begin{equation}
\begin{split}
\F_{211}\F_{131} &= \F_{211464} + \F_{214164} + \F_{214614} + \F_{214641} +
\F_{241164} + \F_{241614} \\
&+ \F_{241641} + \F_{246114} + \F_{246141} + \F_{246411} + \F_{421164} +
\F_{421614} + \F_{421641} \\
&+ \F_{426114} + \F_{426141} + \F_{426411} + \F_{462114} + \F_{462141} +
\F_{462411} + \F_{464211}\,. \\
\end{split}
\end{equation}

\bigskip
Recall that the prime parking functions are those that do not occur
in the decomposition of any nontrivial product $\F_{\park'}\F_{\park''}$.

\subsection{The coalgebra $\PQSym$}

There is a coproduct on $\PQSym$ which appears as a natural extension of the
coproduct of $\FQSym$. Recall (see~\cite{MR,NCSF6}) that if $\sigma$ is
a permutation,
\begin{equation}
\label{coprodF}
\Delta\F_{\sigma} = \sum_{u\cdot v=\sigma}{\F_{\Std(u)} \otimes \F_{\Std(v)}},
\end{equation}
where $\Std$ denotes the usual notion of standardization of a word.

Given a word $w$ on $\{1,2,\ldots\}$, it is possible to define a notion of
\emph{parkization} $\Park(w)$, a parking function which coincides with
$\Std(w)$ when $w$ is a word without repeated letters.

\begin{algorithm}
\label{parkisation}
~

\noindent
\emph{Input}:  A word $w$.

\noindent
\emph{Output}:  A parking function.

Let $n$ be the length of $w$.
Define
\begin{equation}
\label{dw}
d(w):=\min \{i\ |\ \ \  |\{w_j\leq i\}|<i \}\,.
\end{equation}

\begin{itemize}
\item If $d(w)=n+1$, return $w$.
\item Otherwise, let $w'$ be the word obtained by decrementing all
the elements of $w$ greater than $d(w)$. Then return the parkized word
of $w'$.
\end{itemize}
\end{algorithm}

The algorithm is correct since $d(w)=n+1$ iff $w$ is a parking function and
since $w'$ is smaller than $w$ in the lexicographic order, it terminates.

\smallskip
For example, the following tableau displays an execution of the parkization
algorithm: on each line, there is a word $w$ and the value of $d(w)$ and the
next line contains the element $w'$ as defined in the algorithm.

\begin{equation*}
\label{ex-park}
\begin{array}{|ccccccccc||c|}
\hline
&&&&w&&&&& d(w) \\
\hline
\hline
5&7&3&3&13&1&10&10&4 & 2  \\
\hline
4&6&2&2&12&1&9 &9 &3   & 7  \\
\hline
4&6&2&2&11&1&8 &8 &3   & 7  \\
\hline
4&6&2&2&10&1&7 &7 &3   & 9  \\
\hline
4&6&2&2&9& 1&7 &7 &3    & 10 \\
\hline
\end{array}
\end{equation*}

\bigskip
We can now define a coproduct on $\PQSym$ by
\begin{equation}
\label{coprod}
\Delta \F_{\park}:= \sum_{u\cdot v=\park} \F_{\Park(u)} \otimes \F_{\Park(v)}.
\end{equation}

For example,
\begin{equation}
\Delta\F_{121} = 1\otimes\F_{121} + \F_{1}\otimes\F_{21} +
\F_{12}\otimes\F_{1} + \F_{121}\otimes1\,.
\end{equation}
\begin{equation}
\Delta\F_{131} = 1\otimes\F_{131} + \F_{1}\otimes\F_{21} +
\F_{12}\otimes\F_{1} + \F_{131}\otimes1\,.
\end{equation}
\begin{equation}
\Delta\F_{3132} = 1\otimes\F_{3132} + \F_{1}\otimes\F_{132} +
\F_{21}\otimes\F_{21} + \F_{212}\otimes\F_{1} + \F_{3132}\otimes 1\,.
\end{equation}
\begin{equation}
\begin{split}
\Delta\F_{1643165} &= 1\otimes\F_{1643165} + \F_{1}\otimes\F_{532154}
+ \F_{12}\otimes\F_{32154} + \F_{132}\otimes\F_{2143} \\
&+ \F_{1432}\otimes\F_{132} + \F_{15431}\otimes\F_{21}
+ \F_{154315}\otimes\F_{1} + \F_{1643165}\otimes1\,.
\end{split}
\end{equation}

\begin{proposition}
The operation defined by Equation~(\ref{coprod}) is coassociative and is a
morphism for the product. So $(\PQSym,\cdot ,\Delta)$ is a bialgebra.
\end{proposition}

\begin{proof}
The operation is obviously coassociative since the deconcatenation is
coassociative.
Consider two words $w_1$ and $w_2$ and a prefix $u_1$ (resp. $u_2$) of $w_1$
(of $w_2$). Then the set of the parkized words of all prefixes of $w_1\ssh w_2$
containing only letters of $u_1$ and $u_2$ is equal to 
$u_1\ssh u_2$.
So $\Delta$ is a morphism for the product, and hence $\PQSym$ is a bialgebra.
\end{proof}

\subsection{The Hopf algebra $\PQSym$}

Since $\PQSym$ is endowed with a bialgebra structure naturally graded by the
length of parking functions, one defines the antipode as the inverse of the
identity for the convolution product and then endow $\PQSym$ with a Hopf
algebra structure.

The standard formula for the antipode, written on the basis $(\F_{\park})$
reads as
\begin{equation}
\nu (\F_{\park}) = \sum_{r ; u_1\cdots u_r=\park ; |u_i|\geq1} (-1)^r\,
\F_{\Park(u_1)}\F_{\Park(u_2)}\cdots \F_{\Park(u_r)}
\end{equation}

For example,
\begin{equation}
\nu (\F_{122}) = -\F_{122} + \F_{1}\F_{11} + \F_{12}\F_{1} - \F_{1}^3 =
\F_{212} + \F_{221} - \F_{213} - \F_{231} - \F_{321}\,.
\end{equation}

\subsection{The graded dual $\PQSym^*$}

Let $\G_{\park}=\F_{\park}^* \in\PQSym^*$ be the dual basis of $(\F_\park)$.

\begin{proposition}
The product on $\PQSym^*$ is given by
\begin{equation}
\label{prodG}
\G_{\park'} \G_{\park''} = \sum_{\park \in \park'\convP\park''} \G_\park\,,
\end{equation}
where the \emph{convolution} $\park'\convP\park''$ of two parking functions
is defined as
\begin{equation}
\park'\convP\park'' = \sum_{u,v ;
\park=u\cdot v\,\in\,\PF, \Park(u)=\park', \Park(v)=\park''} \park\,.
\end{equation}
\end{proposition}

\begin{proof}
If $\langle\,,\,\rangle$ denotes the duality bracket, the product on
$\PQSym^*$ is given by
\begin{equation}
\G_{\park'} \G_{\park''}
= \sum_{\park}
    \langle\, \G_{\park'}\otimes\G_{\park''}, \Delta\F_\park \,\rangle\,
    \G_\park \\
= \sum_{\park\in\park'\convP\,\park''} \G_\park.
\end{equation}
\end{proof}

For example,
\begin{equation}
\G_{1}\G_{1} = \G_{11} + \G_{12} + \G_{21}, \qquad
\G_{1}\G_{11} = \G_{111} + \G_{122} + \G_{211} + \G_{311}\,.
\end{equation}
\begin{equation}
\G_{1}\G_{12} = \G_{112} + \G_{113} + \G_{123} + \G_{212} + \G_{213}
+\G_{312}\,.
\end{equation}
\begin{equation}
\G_{1}\G_{21} = \G_{121} + \G_{131} + \G_{132} + \G_{221} + \G_{231}
+\G_{321}\,.
\end{equation}
\begin{equation}
\begin{split}
\G_{12} \G_{11} &= \G_{1211} + \G_{1222} + \G_{1233} + \G_{1311} + \G_{1322}\\
&+ \G_{1411} + \G_{1422} + \G_{2311} + \G_{2411} + \G_{3411}\,.
\end{split}
\end{equation}
\begin{equation}
\begin{split}
\G_{211} \G_{131} &=
\G_{211131} + \G_{211141} + \G_{211151} + \G_{211161} + \G_{211242} \\
&+
\G_{211252}
+ \G_{211262} + \G_{211353} + \G_{211363} + \G_{211464} + \G_{322131}
+ \G_{322141} \\
&+
\G_{322151} + \G_{322161} + \G_{433141} + \G_{433151} + \G_{433161} +
\G_{433131} + \G_{544131}\,.
\end{split}
\end{equation}

When restricted to permutations, the product of $\G$ coincides with the
convolution of~\cite{Re,MR}.
Notice also that
\begin{equation}
\G_1^n = \sum_{\park\in\PF_n} \G_\park\,.
\end{equation}

\begin{proposition}
The coproduct $\Delta\G_\park$ is given by
\begin{equation}
\label{coprodG}
\Delta \G_\park := \sum_{u,v ; \park\in u\ssh v}
                   {\G_{u} \otimes \G_{v}}\,.
\end{equation}
\end{proposition}

\begin{proof}
If $\langle\,,\,\rangle$ denotes the duality bracket, the coproduct on
$\PQSym^*$ is given by
\begin{equation}
\Delta\G_{\park} =
 \sum_{\park',\park''}\langle\,\G_{\park}, \F_{\park'} \F_{\park''}\,\rangle\,
    \G_{\park'}\otimes\G_{\park''} 
= \sum_{\park\in \park'\ssh\park''}\G_{\park'}\otimes\G_{\park''}.
\end{equation}
\end{proof}

For example,
\begin{equation}
\label{deltaG121}
\Delta\G_{121} = 1\otimes\G_{121} + \G_{121}\otimes1\,.
\end{equation}
\begin{equation}
\label{deltaG131}
\Delta\G_{131} = 1\otimes\G_{131} + \G_{11}\otimes\G_{1} + \G_{131}\otimes1\,.
\end{equation}
\begin{equation}
\Delta\G_{3132} = 1\otimes\G_{3132} + \G_{1}\otimes\G_{221} +
\G_{12}\otimes\G_{11} + \G_{3132}\otimes 1\,.
\end{equation}
\begin{equation}
\label{164821657}
\begin{split}
\Delta\G_{164821657} =&\ 1\otimes\G_{164821657}
 + \G_{121}\otimes\G_{315324} \\
&+ \G_{1421}\otimes\G_{24213} + \G_{14215}\otimes\G_{1312}
+ \G_{164821657}\otimes1\,.
\end{split}
\end{equation}

There is also a direct way to describe the coproduct of $\G_\park$
in terms of breakpoints:

\begin{proposition}
Let $\park$ be a parking function of length $n$.
For $b$ in $\{0,\ldots,n\}$, define $\park'(b)$ and $\park''(b)$ as the
restrictions of $\park$ to the respective intervals $[1,b]$ and $[b+1,n]$.
Then
\begin{equation}
\Delta \G_\park := \sum_{b}
                   {\G_{\park'(b)} \otimes \G_{\park''(b)}}\,,
\end{equation}
where the sum runs over all breakpoints of $\park$.
\end{proposition}

\begin{proof}
The term $\G_u\otimes\G_v$ appears in $\Delta\G_\park$ iff $\park$
belongs to the shifted shuffle of $u$ and $v$, so that $\park$ has a
breakpoint at $b=|u|$. Then $\a'(b)=u$ and $\a''(b)=v$.
\end{proof}

For example, the breakpoints of $164821657$ are $\{0,3,4,5,9\}$ so that one
recovers the result of Equation~(\ref{164821657}).

Let us finally mention that the realization provided in
Section~\ref{realG-sec} allows to shed an interesting light on the coproduct
and the fact that $\PQSym^*$ is a Hopf algebra.

\subsection{$\PQSym$ as a bidendriform bialgebra}

In~\cite{Foi}, Foissy has proved that the Hopf algebra of Free quasi-symmetric
functions $\FQSym$ is bidendriform. A very slight modification of his
operations allows us to state:

\begin{theorem}
$\PQSym^*$ is a bidendriform bialgebra with the following definitions:
\begin{equation}
\G_{\park'} \gaudend \G_{\park''} =
\sum_{\park=u.v\in \park'\convP\park'', |u|=|\park'| ; \max(v)<\max(u)}
\G_\park,
\end{equation}
\begin{equation}
\G_{\park'} \droitdend \G_{\park''} =
\sum_{\park=u.v\in \park'\convP\park'', |u|=|\park'| ; \max(v)\geq\max(u)}
\G_\park,
\end{equation}
\begin{equation}
\Delta_\gaudend \G_\park = \sum_{\park\in u\ssh v ; \last(\park)\leq |u|}
\G_u\otimes \G_v,
\end{equation}
\begin{equation}
\Delta_\droitdend \G_\park = \sum_{\park\in u\ssh v ; \last(\park)>|u|}
\G_u\otimes \G_v.
\end{equation}
where $|u|\geq1$ and $|v|\geq1$, and $\last(\a)$ means the
last letter of $\a$.
\end{theorem}

\begin{proof}
First, the three defining relations of a dendriform dialgebra are satisfied.
Let us check the first one, for instance. The left part of (\ref{dend1})
amounts to consider the elements $w$ in $\a\convP\b\convP\c$ where the last
maximum of $w$ belongs to $\a$. It is the same for the right part of
(\ref{dend1}). The other two relations are proved in the same way, by checking
that they build the words in $\a\convP\b\convP\c$ where the last maximum is in
$\b$ (Equation~\ref{dend2}) or in $\c$ (Equation~\ref{dend3}). 

The three defining relations of a codendriform coalgebra are also satisfied
since they amount to split the set of parking functions indexing the
elements of $(\overline\Delta\otimes Id)\circ\overline\Delta(\G_\park)$
according to the element of the tensor product containing the last letter of
$\park$.

Since the sum of the four compatibility relations is equivalent to the
coassociativity of $\overline\Delta$, it is sufficient to check any three of
them. We will only prove the first one (more complicated than the second and
third one) in detail, the other ones being proved in the same way.

We will identify till the end of this proof any function $\G_\park$ with its
index $\park$.
Let $\park$ and $\parkb$ be two parking functions of length $p$ and $q$.
For any word $w$ of length $p+q$, let $m_1=\min(w_1\ldots w_{p})$,
$m_2=\min(w_{p+1}\ldots w_{p+q})$, $M_1=\max(w_1\ldots w_{p})$,
$M_2=\max(w_{p+1}\ldots w_{p+q})$.

Let $i$ be any integer and split into two groups the parking functions
indexing the terms in $\G_{\park}\droitdend\G_{\parkb}$ having a breakpoint at
$i$ according to the criterions:
\begin{itemize}
\item $(M_1\leq i)$, or $(M_1> i$  and $m_1\leq i)$, or $(m_1>i)$,
\item $(m_2\leq i)$, or $(m_2> i)$.
\end{itemize}

Apply $\Delta_\droitdend$ to the first group and consider $S_1$ as the sum of
the elements $\G_u\otimes\G_v$ such that $u$ is of length $i$. Since $M_1\leq i$
and $m_2\leq i$, the first $p$ letters of all words are in $u$, whereas the
others are both in $u$ and $v$ (by hypothesis, the last one is in $v$).  Note
that the positions of the letters belonging to the right-hand side of the
tensor product are independent of the element of the first group and are the
positions of the $p+q-i$ greatest letters of $\parkb$.
Moreover, there exists a breakpoint of $\parkb$ separating those letters from
the other ones. Since the last letter of any word of the first group goes to
the right-hand side of the tensor product and comes from the last letter of
$\parkb$, we then deduce that the right-hand side of $S_1$ is built with all
$\parkb''_\droitdend$ of length $p+q-i$, using Sweedler's notations.
Finally, in the left-hand side, we find the elements of the form
$\park\ssh\parkb'_\droitdend$, and all those ones, since they
correspond to the restriction of all words of $\park\droitdend\parkb$ to
letters smaller than $i$. Finally, summing up over all possible $i$, we get
that the sum of all the elements of all the first groups is
$\G_{\park}\G_{\parkb'_\droitdend}\otimes \G_{\parkb''_\droitdend}$,
that is, the fourth term of the right hand-side of Equation~(\ref{bidend1}).

In the same way, one proves that the second group, corresponding to
$M_1\leq i$ and $m_2>i$ gives the term $\park\otimes\parkb$ of
Equation~(\ref{bidend1}).
The third group, corresponding to $M_1>i$, $m_1\leq i$, and $m_2\leq i$ gives
the term
$\park'\parkb'_\droitdend\otimes \park''\droitdend\parkb''_\droitdend$ of 
Equation~(\ref{bidend1}).
The fourth group, corresponding to $M_1>i$, $m_1\leq i$, and $m_2>i$ gives
the term $\park'\otimes \park''\droitdend\parkb$ of
Equation~(\ref{bidend1}).
The fifth group, corresponding to $m_1>i$ and $m_2\leq i$ gives the term
$\parkb'_\droitdend \otimes \park\droitdend\parkb''_\droitdend$ of
Equation~(\ref{bidend1}).
The sixth group, corresponding to $m_1>i$ and $m_2>i$, gives no term since we
would have $|u|=0$, which is impossible.
\end{proof}

For example,
\begin{equation}
\label{pqs-gauche}
\begin{split}
\G_{12}\gaudend\G_{212} =& \ \ \ \
 \G_{13212} + \G_{14212} + \G_{14313} + \G_{14323} + \G_{15212} +
 \G_{15313} \\ &+
 \G_{15323} + \G_{24313} + \G_{24212} + \G_{34212} + \G_{23212} +
 \G_{25212} \\ &+
 \G_{25313} + \G_{35212} + \G_{45212}.
\end{split}
\end{equation}
\begin{equation}
\label{pqs-droite}
\begin{split}
\G_{12}\droitdend\G_{212} = & \ \ \ \
 \G_{12212} + \G_{12313} + \G_{12323} + \G_{12414} + \G_{12424} +
 \G_{12434} \\ &+
 \G_{13313} + \G_{13323} + \G_{13414} + \G_{13424} + \G_{23313} + \G_{23414}.
\end{split}
\end{equation}
\begin{equation}
\Delta_\gaudend \G_{1252754} = \G_{125254}\otimes\G_1
 + \G_{1224}\otimes\G_{131}.
\end{equation}
\begin{equation}
\Delta_\droitdend \G_{1252754} = \G_{122}\otimes \G_{2421}
 + \G_1\otimes\G_{141643}.
\end{equation}

The duality of bidendriform bialgebras implies that the bidendriform
relations for $\PQSym$ are
\begin{equation}
\F_{\park'} \gaudend \F_{\park''} =
\sum_{\a\in \a'\ssh \a'' ; \last(\a)\leq|\a'|} \F_\park,
\end{equation}
\begin{equation}
\F_{\park'} \droitdend \F_{\park''} =
\sum_{\a\in \a'\ssh \a'' ; \last(\a)>|\a'|}
\F_\park,
\end{equation}
\begin{equation}
\Delta_\gaudend  \F_\park = \sum_{u\cdot v=\a ; \max(v)<\max(u)}
\F_{\Park(u)}\otimes \F_{\Park(v)},
\end{equation}
\begin{equation}
\Delta_\droitdend  \F_\park = \sum_{u\cdot v=\a ; \max(v)\geq\max(u)}
\F_{\Park(u)}\otimes \F_{\Park(v)},
\end{equation}
where the sums inside the coproducts occur over non trivial deconcatenations,
that is $|u|\geq1$ and $|v|\geq1$.

\smallskip
We then have the following consequences of the results of Foissy~\cite{Foi}. 

\begin{corollary}
$\PQSym$ is a self-dual Hopf algebra.
\end{corollary}

\begin{corollary}
The Lie algebra of primitive elements of $\PQSym$ is a free Lie algebra.
\end{corollary}

Let
\begin{equation}
PF(t)=1+\sum_{n\ge 1}(n+1)^{n-1}t^n\,.
\end{equation}

\begin{corollary}
$\PQSym$ is free as a dendriform dialgebra on its totally primitive elements
whose degree generating series is
\begin{equation}
\begin{split}
TP(t) :=& \frac{(PF(t)-1)}{PF(t)^2} \\
=&
t + t^2 + 7\, t^3 + 66\, t^4 + 786\, t^5 + 11\,278\, t^6 + 189\,391\, t^7 \\
& + 3\,648\,711\, t^8 + 79\,447\,316\, t^9 + O(t^{10}).
\end{split}
\end{equation}
\end{corollary}

For example, $\F_{1}$ and $\G_{1}$ are totally primitive and so are
$\F_{12}-\F_{11}$ and $\G_{11}$.
Here are bases of the seven dimensional space of totally primitive
elements of $\PQSym$ and $\PQSym^*$ in degree $3$:

\begin{equation}
\begin{split}
&
\F_{123}-\F_{122}-\F_{112}+\F_{111},\quad
 \F_{311}-\F_{211}, \quad
 \F_{113}-\F_{112}, \quad \\ &
 \F_{131}-\F_{121}, \quad
 \F_{132}-\F_{131}, \quad
 \F_{231}-\F_{121}, \quad
 \F_{213}-\F_{212}.
\end{split}
\end{equation}
\begin{equation}
\begin{split}
&\G_{122}-\G_{212},\quad
 \G_{131}-\G_{311}, \quad
 \G_{312}-\G_{132}, \quad \\ &
 \G_{111}, \quad
 \G_{112}, \quad
 \G_{121}, \quad
 \G_{211}.
\end{split}
\end{equation}

Thanks to the bidendriform structure of $\PQSym$, we know that $\PQSym$ and
$\PQSym^*$ are isomorphic as bidendriform bialgebras and hence isomorphic as
Hopf algebras. We do not know an explicit isomorphism, but restricting to
$\FQSym$, that is, permutations,  the linear map $\varphi$ defined by
\begin{equation}
\label{phi-isom}
\varphi ( \F_\sigma) := \sum_{\park, \Std(\park)=\sigma^{-1}}{\G_\park}\,,
\end{equation}
is a bidendriform and hence a Hopf embedding, compatible with the usual
realization of $\FQSym$~\cite{NCSF6}.

\subsection{Free generators and multiplicative bases}
\label{alg-sec}

Let us say that a word $w$ over $\NN^*$ is \emph{connected} if it cannot be
written as a shifted concatenation $w=u\sconc v$, and \emph{anti-connected} if
its mirror image $\overline{w}$ is connected.

\begin{proposition}
\label{PQS-free}
$\PQSym$ is free over the set
\begin{equation}
\left\{ \F_{\bf c}\, |\, {\bf c}\in\PF, \c \text{\ connected} \right\}
\end{equation}
and $\PQSym^*$ is free over the set
\begin{equation}
\left\{ \G_{\bf d}\, |\, {\bf d}\in\PF, \d \text{\ anti-connected} \right\}
\end{equation}
\end{proposition}

\begin{proof}
Clearly, any word $w$ has a unique \emph{maximal factorization} into
connected words, $w = w_1\sconc w_2 \sconc\ldots\sconc w_k$ where all $w_i$
are connected.
Moreover, the lexicographically minimal word in $w_1\ssh\ldots\ssh w_k$ is $w$
so that the matrix expressing all products of $\F$ indexed by connected
words is triangular over the basis $\F_\a$, with ones on the diagonal. 
The proof is exactly the same for the $\G$.
\end{proof}

The ordinary generating function for the numbers $c_n$ of connected parking
functions is
\begin{equation} 
  \begin{split}
    \sum_{n\geq1} c_n t^n &= 1 - PF(t)^{-1} \\
    &= t + 2\, t^2 + 11\, t^3 + 92\, t^4 + 1\,014\, t^5 + 13\,795\, t^6
       + 223\,061\, t^7 \\ &
       + 4\,180\,785\, t^8 
      + 89\,191\,196\, t^9  + 2\,135\,610\,879\, t^{10}
      + 56\,749\,806\,356\, t^{11}\\ &+ 1\,658\,094\,051\,392\, t^{12}
      + O\left ({t}^{13}\right )\,.
  \end{split}
\end{equation}

Let $\park = {\bf a}_1 \sconc {\bf a}_2\sconc \cdots\sconc {\bf a}_r$ be the
maximal factorization of $\park$ into connected parking functions. We set
\begin{equation}
\F^{\park} = \F_{\park_1} \cdot \F_{\park_2} \cdots \F_{\park_r}\,,
\end{equation}
and
\begin{equation}
\G^{\overline{\park}} = \G_{\overline{\park_r}} \cdots
 \G_{\overline{\park_1}}\,.
\end{equation}

\begin{proposition}
The basis $(\F^{\park})$ of $\PQSym$ and the basis $(\G^{\overline{\park}})$
of $\PQSym^*$ are both multiplicative.
\end{proposition}

\begin{proof}
This follows from the proof of Proposition~\ref{PQS-free}.
\end{proof}

\medskip
Now, if $\S_\park$ (resp. $\T_\park$) is the dual basis of $\F^{\park}$
(resp. $\G^{\overline{\park}}$) then
\begin{equation}
\{ \S_{\bf c} \,|\, {\bf c} \text{\ connected} \} \text{\ and\ }
\{ \T_{\bf c} \,|\, {\bf c} \text{\ connected} \}
\end{equation}
are bases of the primitive Lie algebras $\LPQ^*$ (resp. $\LPQ$)
of $\PQSym^*$ (resp. $\PQSym$). 

Thanks again to \cite{Foi}, we known that both Lie algebras are free, on
generators whose degree generating series is
\begin{equation}
  \begin{split}
    1 - \prod_{n\geq1}{(1-t^n)}^{c_n} &= 1-(1-t)(1-t^2)^2(1-t^3)^{11} \cdots\\
    &= t + 2\,t^2 + 9\,t^3 + 80\,t^4 + 901\,t^5 + 12\,564\,t^6
      + 206\,476\,t^7 \\
    & + 3\,918\,025\,t^8 + 84\,365\,187\,t^9 + 2\,034\,559\,143\,t^{10} +
      O\left ({t}^{11}\right )\,.
  \end{split}
\end{equation}

\subsection{$\PQSym^*$ as a combinatorial Hopf algebra}

Since $\FQSym$ can be embedded in $\PQSym$, we have a canonical Hopf embedding
of $\NCSF$ in $\PQSym$ given by
\begin{equation}
S_n \mapsto \F_{12\cdots n}\,.
\end{equation}
With parking functions, we have other possibilities: for example,
\begin{equation}
j(S_n) := \F_{11\cdots 1}
\end{equation}
is a Hopf embedding, whose dual $j^*$ maps $\PQSym^*$ to $\QSym$ and
therefore endows $\PQSym^*$ with a different structure of combinatorial Hopf
algebra in the sense of~\cite{ABS}.

On the dual side, we have a Hopf embedding
\begin{equation}
S_n \mapsto \sum_{\Std(\park)=12\cdots n} \G_\park
\end{equation}
of $\NCSF$ into $\PQSym^*$, given by  the restriction of the
self-duality isomorphism of Formula~(\ref{phi-isom}) to the $\NCSF$ subalgebra
$S_n=\F_{12\cdots n}$ of $\PQSym$. Its transpose gives a Hopf epimorphism
$\eta:\ \PQSym\rightarrow QSym$, which maps $\F_\a$ to $F_I$, where
$I$ is the descent composition of the word $\a$.

\pagebreak[4]
\subsection{Realizations of $\PQSym^*$ and $\PQSym$}

\subsubsection{Realization of $\PQSym^*$}
\label{realG-sec}

The algebra $\PQSym^*(A)$ admits a simple realization in terms of
noncommutative polynomials \cite{NT2}, which is similar to the construction of
$\FQSym$.
If $A$ is a totally ordered infinite alphabet, one can set
\begin{equation}
\label{realG}
\G_\park(A) :=\sum_{w\in A^*, \Park(w)=\park}w.
\end{equation}

\begin{theorem}[\cite{NT2}]
\label{realG-thm}
These polynomials satisfy Relations (\ref{prodG}) and allow to write the
coproduct as $\Delta \G_\park = \G_\park(A'\hat{+} A'')$ where $A'\hat{+} A''$
denotes the ordered sum of two mutually commuting alphabets isomorphic to $A$
as ordered sets.
\end{theorem}

Let us recall the precise way to introduce a coproduct on an algebra realized
on words under certain conditions. Start with $A'$ and $A''$, two mutually
commuting alphabets isomorphic to $A$ as ordered sets. Then build their
ordered sum $A'\hat{+} A''$ and compute $\G_\park(A'\hat{+} A'')$ separating
inside each term what belongs to $A'$ and what belongs to $A''$.
Assume that one can write, for all parking function $\park$,
\begin{equation}
\G_\park(A'\hat{+} A'') = \sum_{\park',\park''} \G_{\park'}(A')\G_{\park''}(A'')
\end{equation}
where the sum is taken over a set of pairs of parking functions depending on
$\park$.
Then the operation
\begin{equation}
\Delta\G_{\park} := \sum_{\park',\park''} \G_{\park'}\otimes\G_{\park''},
\end{equation}
where the sum is taken over the same set as before is a coproduct.

For example, $\G_{121}=\sum_{i} a_ia_{i+1}a_i$, so that
\begin{equation}
\G_{121}(A'\hat{+} A'') = \G_{121}(A') + \G_{121}(A''),
\end{equation}
since $a_i\in A'$ is equivalent to $a_{i+1}\in A'$ by definition of the
ordered sum of alphabets. One then recovers the results of
Equation~(\ref{deltaG121})
Now, $\G_{131}=\sum_{i,j ; j>i+1}a_ia_ja_i$, so that
\begin{equation}
\G_{131}(A'\hat{+} A'') = \G_{131}(A') + \G_{11}(A')\G_{1}(A'') +
\G_{121}(A''),
\end{equation}
since $a_i$ and $a_j$ can belong to $A$ (first term), or $a_i$ belongs to $A'$
and $a_j$ belongs to $A''$ (second term), or $a_i$ and $a_j$ belong to $A''$
(third term). One then recovers the results of Equation~(\ref{deltaG131}).

\subsubsection{Realization of $\PQSym$}

Although $\PQSym$ and $\PQSym^*$ are isomorphic as Hopf algebras, 
no explicit isomorphism is known.   We can nevertheless 
propose  a realization of
$\PQSym$ in terms of $(0,1)$-matrices instead of words.

This construction is reminiscent of the construction of $\MQSym$
(see~\cite{Hiv,NCSF6}), and coincides with it when restricted to
permutation matrices, providing the natural embedding of $\FQSym$ in
$\MQSym$.

Let $\MA_n$ be the vector space spanned by symbols $X_M$ where $M$ runs over
$(0,1)$-matrices with $n$ columns and an infinite number of rows, with $n$
nonzero entries, so that at most $n$ rows are nonzero.

Given such a matrix $M$, we define its \emph{vertical packing}
$\vp(M)$ as the finite matrix obtained by removing the null rows of $M$.

For a vertically packed matrix $P$, we define
\begin{equation}
\M_P = \sum_{\vp(M)=P} X_M\,.
\end{equation}

Now, given a $(0,1)$-matrix, we define its reading $r(M)$ as the word
obtained by reading its entries by rows, from left to right and top to bottom
and recording the numbers of the columns of the ones. For example, the reading
of the matrix
\begin{equation}
\begin{pmatrix}
0 & 1 & 1 & 0 \\
1 & 0 & 0 & 0 \\
0 & 1 & 0 & 0
\end{pmatrix}
\end{equation}
is $(2,3,1,2)$.

A matrix $M$ is said to be of \emph{parking type} if $r(M)$ is a parking
function.
Finally, for a parking function $\park$, we set
\begin{equation}
\F_{\park} := \sum_{r(P)=\park, \text{$P$ vertically packed}}{\M_P} =
\sum_{r(M)=\park} X_M\,.
\end{equation}

For example,
\begin{equation}
\F_{(1,2,2)} = \M_{
\begin{pmatrix}
1 & 1 & 0 \\
0 & 1 & 0 \\
\end{pmatrix}
}  + \M_{
\begin{pmatrix}
1 & 0 & 0 \\
0 & 1 & 0 \\
0 & 1 & 0 \\
\end{pmatrix}
}\,.
\end{equation}

The multiplication on $\MA=\bigoplus_{n} \MA_n$ is defined by columnwise
concatenation of the matrices:
\begin{equation}
X_M X_N = X_{M\cdot N}\,.
\end{equation}

In order to explicit the product of $\M_P$ by $\M_Q$, we first need a
definition.
Let $P$ and $Q$ be two vertically packed matrices with respective heights $p$
and $q$. The \emph{augmented shuffle} of $P$ and $Q$ is defined as follows:
let $r$ be an integer in $[\max(p,q),p+q]$. One inserts zero rows in $P$ and
$Q$ in all possible ways so that the resulting matrices have $p+q$ rows.
Let $R$ be the matrix obtained by concatenation of such pairs of matrices. The
augmented shuffle consists in the set of such matrices $R$ with nonzero rows.
We denote this set by $\saug(P,Q)$.

\begin{theorem}
The following formulas hold:
\begin{equation}
\label{matMprod}
\M_P \M_Q = \sum_{R\,\in\,\saug(P,Q)} {\M_R}\,,
\end{equation}
and
\begin{equation}
\label{realFprod}
\F_{\park'}\F_{\park''}=\sum_{\park\in\park'\,\ssh\,\park''}\F_\park\,,
\end{equation}
\end{theorem}
This  is the same as Equation~(\ref{prodF}).
\begin{proof}
Formula (\ref{matMprod}) comes from the definition of the augmented shuffle of
matrices: any matrix in $\saug(P,Q)$ appears as a product $X_M.X_N$ where
$vp(M)=P$ and $vp(N)=Q$.
Conversely, any element in $\M_P\M_Q$ has as vertical packing a matrix with a
number of rows in the interval $[\max(p,q),p+q]$ which left part has $M$ as
vertical packing and right part has $N$ as vertical packing.

The proof of (\ref{realFprod}) is almost the same as the previous one
if one starts from the definition $\F_\park = \sum_{r(M)=\park} X_M$.
\end{proof}

\smallskip
Finally, concerning the coproduct, one has first to define the parkization
$\Park(M)$ of a vertically packed matrix $M$, which consists in iteratively
removing column $d(r(M))$ until $M$ becomes a parking matrix.

The coproduct of a matrix $\M_P$ is then defined as:
\begin{equation}
\Delta \M_{P} = \sum_{Q\cdot R=P} \M_{\Park(Q)} \otimes \M_{\Park(R)}\,,
\end{equation}
It is then easy to check that

\begin{proposition}
The following formula holds:
\begin{equation}
\Delta \F_{\park} =
  \sum_{u\cdot v=\park} \F_{\Park(u)} \otimes \F_{\Park(v)}\,.
\end{equation}
\end{proposition}
This is the same as Equation~(\ref{coprodF}).

\subsubsection{Realization of $\FQSym$}

A parking matrix $M$ is said to be a \emph{word matrix} if there is exactly
one $1$ in each column.
Then $\FQSym$ is the Hopf subalgebra generated by the parking word matrices.

\subsection{$\PQSym^*$ as a dendriform trialgebra}

Since we already know that $\K\langle A\rangle^+$ is a dendriform trialgebra
(see Definition~\ref{def-trigdend} and Lemma~\ref{lem-trigdend}), and since
$\PQSym^*$ can be realized on words, it is a natural question to ask whether
$\PQSym^*$ is a sub-trialgebra of $\K\langle A\rangle^+$.

\begin{theorem}
\label{th-pqstrig}
$\PQSym^*$ is a sub-dendriform trialgebra of $\K\langle A\rangle^+$ with the
following product rules:
\begin{equation}
\G_{\park'} \gautrid \G_{\park''} =
\sum_{\park=u.v\in \park'\convP\park'', |u|=|\park'| ; \max(v)<\max(u)}
\G_\park,
\end{equation}
\begin{equation}
\G_{\park'} \miltrid \G_{\park''} =
\sum_{\park=u.v\in \park'\convP\park'', |u|=|\park'| ; \max(v)=\max(u)}
\G_\park,
\end{equation}
\begin{equation}
\G_{\park'} \droittrid \G_{\park''} =
\sum_{\park=u.v\in \park'\convP\park'', |u|=|\park'| ; \max(v)>\max(u)}
\G_\park,
\end{equation}
\end{theorem}

\Proof
Since $\PQSym^*$ can be realized on words, one only needs to check that
$\PQSym^*$ is stable under all three operations, their compatibility coming
from the fact that $\K\langle A\rangle^+$ is a dendriform trialgebra.
Since all words having a given parkized word have the same inversions, they
have in particular the same relations between the maximum of any prefix and
any suffix of given lengths.
One then derives the product rules from direct calculation.
\qed

For example, with the notations of Definition~\ref{def-trigdend}:
\begin{equation}
\begin{split}
\G_{12}\gautrid\G_{212} =& \ \ \ \
 \G_{13212} + \G_{14212} + \G_{14313} + \G_{14323} + \G_{15212} +
 \G_{15313} \\ &+
 \G_{15323} + \G_{24313} + \G_{24212} + \G_{34212} + \G_{23212} +
 \G_{25212} \\ &+
 \G_{25313} + \G_{35212} + \G_{45212}.
\end{split}
\end{equation}
\begin{equation}
\G_{12}\miltrid \G_{212} = \G_{12212} + \G_{13313} + \G_{13323} + \G_{23313}.
\end{equation}
\begin{equation}
\begin{split}
\G_{12}\droittrid\G_{212} = & \ \ \ \
 \G_{12313} + \G_{12323} + \G_{12414} + \G_{12424} \\ &+
 \G_{12434} +
 \G_{13414} + \G_{13424} + \G_{23414}.
\end{split}
\end{equation}

Based on numerical evidence, we conjecture the following result:

\begin{conjecture}
$\PQSym^*$ is a free dendriform trialgebra.
\end{conjecture}

Recall that the generating series $F(t)$ for the dimensions of the free
dendriform trialgebra satisfies
\begin{equation}
F(t)-1=t (2F(t)^2-F(t)).
\end{equation}
Applying the same trick as in~\cite{Foi} for computing the generating
series of the totally primitive elements, one gets the generating series of
the number $g_n$ of generators in degree $n$ of $\PQSym^*$ as a free
dendriform trialgebra:

\begin{equation}
\label{2gens}
\begin{split}
\sum_{n\geq0} g_nt^n = \frac{PF(t)-1}{2PF(t)^2 - PF(t)} \\
=& t + 5\,t^3 + 50\, t^4 + 634\,t^5 + 9\,475\,t^6 + 163\,843\,t^7 \\
&+ 3\,226\,213\,t^8 + 71\,430\,404\,t^9 + O(t^{10}).
\end{split}
\end{equation}

By self-duality of $\PQSym$, one can endow $\PQSym$ with a structure of
dendriform trialgebra.

Note that $\FQSym$ is not a sub-dendriform trialgebra of
$\K\langle A\rangle^+$ since the product is not internal and that,
independently of the realization, it cannot be a free dendriform trialgebra
since the substitution $F(t)=\sum_{n} {n!t^n}$ in Equation~(\ref{2gens}) does
not yield a series with nonngative integer coefficients.

\subsection{The internal product}

We shall now recall  the definition of the internal product
of $\PQSym$, introduced in \cite{NT2}.
We first need a few standard notations about biwords.
Let $x_{ij}=$ \scriptsize$\binomial{i}{j}$\normalsize\ be commuting
indeterminates, and $a_{ij}=$ \scriptsize$\ncbinomial{i}{j}$\normalsize\
be noncommuting ones. We shall denote by
\scriptsize $\binomial{i_1\ i_2 \cdots i_r}{j_1\ j_2 \cdots j_r}$ \normalsize\
the monomial
\scriptsize $\binomial{i_1}{j_1}\binomial{i_2}{j_2}\cdots \binomial{i_r}{j_r}$
\normalsize\
and by
\scriptsize $\ncbinomial{i_1,i_2,\cdots i_r}{j_1,j_2,\cdots j_r}$ \normalsize\
the word
\scriptsize\
$\ncbinomial{i_1}{j_1}\ncbinomial{i_2}{j_2}\cdots \ncbinomial{i_r}{j_r}$.
\normalsize\
Such expressions will be referred to respectively as \emph{bimonomials} and
\emph{biwords}.

\medskip
Recall that Gessel constructed the descent algebra by extending to $\QSym$ the
coproduct dual to the internal product of symmetric functions. That is, if $X$
and $Y$ are two totally and isomorphically ordered alphabets of commuting
variables, we can identify a tensor product $f\otimes g$ of quasi-symmetric
functions with $f(X)g(Y)$. Denoting by $XY$ the Cartesian product $X\times Y$
endowed with the lexicographic order, Gessel defined for $f\in QSym_n$
\begin{equation}
\delta(f)=f(XY) \in QSym_n\otimes QSym_n \,.
\end{equation}
The dual operation on $\Sym_n$ is the internal product $*$, for which it
is anti-isomorphic to the descent algebra $\Sigma_n$.
This construction can be extended to $\PQSym^*$.
Let $A'$ and $A''$ be two totally and isomorphically ordered alphabets of
noncommuting variables, but such that $A'$ and $A''$ commute with each other.
We denote by $A'A''$ the Cartesian product $A'\!\times\! A''$ endowed with the
lexicographic order. This is a total order in which each element has a
successor, so that $G_\park(A'A'')$ is a well defined polynomial.
Identifying tensor products of words of the same length with
words over $A'A''$, we have
\begin{equation}
\G_\park(A'A'')=\sum_{\Park(u\otimes v)=\park} u\otimes v\,.
\end{equation}
For example, writing tensor products as biwords, one has
\begin{equation}
\G_{4121}(A'A'') = \sum_{a,b,c,d} {\left[\,
\begin{matrix}b&a&a&a\\ d&c&c+1&c\end{matrix}\,\right]}
\end{equation}
with $b>a$, or $b=a$ and $d\geq c+3$.
\begin{theorem}[\cite{NT2}]
The formula $\delta(\G_\park)=\G_\park(A'A'')$ defines a coassociative
coproduct on each homogeneous component $\PQSym_n^*$. Actually,
\begin{equation}
\label{coprodmult}
\delta(\G_\park)=\sum_{\Park(\park'\otimes\park'')=\park}
                     {\G_{\park'}\otimes \G_{\park''}}\,,
\end{equation}
where $\park'$ and $\park''$ are parking functions.
By duality, the formula
\begin{equation}
\F_{\park'} * \F_{\park''} = \F_{\Park(\park'\otimes\park'')}
\end{equation}
defines an associative product on each $\PQSym_n$.
\end{theorem}

Since $A$ is infinite, $\delta$ is compatible with the product of $\PQSym^*$.

\begin{example}
\begin{equation}
\begin{split} 
\delta\G_{4121} =& (\G_{2111}+\G_{3111}+\G_{4111})\otimes
(\G_{1232}+\G_{1121}+\G_{2121}+\G_{3121}+\G_{4121}) \\
&+\G_{1111}\otimes\G_{4121}.
\end{split}
\end{equation}
\end{example}

\begin{example}
\begin{equation}
\F_{211}*\F_{211}=\F_{311}; \qquad \F_{211}*\F_{112} = \F_{312};
\end{equation}
\begin{equation}
\F_{211}*\F_{121}=\F_{321}; \qquad \F_{112}*\F_{312} = \F_{213};
\end{equation}
\begin{equation}
\F_{31143231}*\F_{23571713} = \F_{61385451}.
\end{equation}
\end{example}

Note that although parking functions can be interpreted as chambers
of the Shi arrangement, our internal product is not induced by the
face semigroup of this arrangement. Indeed, one should obtain in
particular an idempotent semigroup, which is clearly not the case.

The main tool for handling internal products of non-commutative symmetric
functions is the splitting formula (see~\cite{NCSF1}, Proposition 5.2). It
does not hold in $\PQSym$, but one can find subalgebras of $\PQSym$ larger
than $\NCSF$ in which it remains true.

\newpage
\section{The Schr\"oder Quasi-Symmetric Hopf algebra $\SQSym$} 
\label{sqsym}

In Section~\ref{trigd1}, we recalled that the little Schr\"oder numbers build
up the Hilbert series of the free dendriform trialgebra on one generator
$\TD$.
We show in~\cite{NTn} that $\TD$ realized on words has a natural structure of
bidendriform bialgebra. In particular, this proves that there is a natural
self-dual Hopf structure on $\TD$.

But parking functions provide another way to find little Schr\"oder numbers.
Indeed, the number of classes of parking functions of length
$n$ under the hypoplactic congruence is also equal to $s_n$. This construction
leads to a non self-dual Hopf algebra, denoted by $\SQSym$.

\subsection{Hypoplactic classes of parking functions}

Let $\congru$ denote the hypoplactic congruence (see~\cite{NCSF4,Nov}).
Recall that the equivalence classes of words under this congruence are
parametrized by \emph{quasi-ribbon tableaux}.
A \emph{quasi-ribbon tableau} of shape $I$ is a ribbon diagram $r$
of shape $I$ filled by letters in such a way that each row of $r$ is
nondecreasing from left to right, and each column of $r$ is strictly
increasing from \emph{top to bottom}.
A word is said to be a \emph{quasi-ribbon word} of shape $I$
if it can be obtained by reading from \emph{bottom to top} and from
left to right the columns of a quasi-ribbon diagram of shape $I$.
For example, the word $11425477$ is a quasi-ribbon word since it is the
reading of the following quasi-ribbon
\begin{equation}
\label{exquasi}
\PetitTableau
\Tableau{1&1&2\\\ &\ &4&4\\\ &\ &\ &5&7&7 \\}
\end{equation}

The hypoplactic classes of parking functions correspond to
\emph{parking quasi-ribbons}, that is, quasi-ribbon words that are parking
functions. We denote this set by $\PQR$, and $\PQR_n$ is the set of
quasi-ribbon parking functions of length $n$.

We will make use of a simple parametrization of the elements of $\PQR$:
define a \emph{segmented word} as a finite sequence of non-empty words,
separated by vertical bars, \emph{e.g.}, $232\sep14\sep5\sep746$.

The parking quasi-ribbons can be represented as segmented nondecreasing
parking functions where the bars only occur at positions
$\cdots a\sep b\cdots$, with $a<b$.
For example, the quasi-ribbon of Equation~(\ref{exquasi}) is represented by
the word $112\sep44\sep577$.

Clearly, a nondecreasing word containing exactly $l$ different letters admits
$2^{l-1}$ segmentations.

On another hand, the statistic $l$ (the length of the packed evaluation
vector) on nondecreasing parking functions has the same distribution as the
number of blocks in non-crossing partitions through the natural bijection.
This is given by a classical $q$-Catalan, $c_n(q)$ (see,
\emph{e.g.},~\cite{Nar}) and finally, the number of canonical packed words of
length $n$ is $c_n(2)$, which is known to be equal to the Schr\"oder number
$s_n$.

For example, $c_1(q)=1$, $c_2(q)=1+q$ and $c_3(q)=1+3q+q^2$, so that
$c_1(2)=1$, $c_2(2)=3$ and $c_3(2)=11$ as one can check on
Equations~(\ref{ex12}) and~(\ref{ex3}).
The coefficients of $c_n(q)$ are known as the Narayana numbers (sequence
A001263 of Sloane's database~\cite{Slo}).

Here is for $n\geq3$ the list of canonical hypoplactic parking functions.
\begin{equation}
\label{ex12}
\{1\}, \qquad\qquad \{11,\, 12,\, 1\sep2\},
\end{equation} 
\begin{equation}
\label{ex3}
\{111\ \ 112\ \ 11\sep2 \ \ 113\ \ 11\sep3\ \ 122\ \ 1\sep22\ \
123\ \ 1\sep23\ \ 12\sep3\ \ 1\sep2\sep3 \}.
\end{equation}

In the sequel, we will identify parking quasi-ribbons and their encodings as
segmented words.

\subsection{The Schr\"oder Quasi-Symmetric Hopf algebra $\SQSym$}

Let us denote by ${\sf P}(w)$ the hypoplactic $P$-symbol of a word $w$ (its
quasi-ribbon). The ${\sf P}$-symbols of parking functions are therefore
parking quasi-ribbons.
With a parking quasi-ribbon $\qrpark$, we associate the elements
\begin{equation}
\PS_{\qrpark} := \sum_{{\sf P}(\park)=\qrpark} {\F_\park},
\quad\text{and}\quad
\QS_{\qrpark} := \overline{\G_\park}\,,
\end{equation}
where $\overline w$ denotes the hypoplactic class of $w$.
For example,
\begin{equation}
\PS_{11|3} = \F_{131}+\F_{311}\,, \qquad
\PS_{113} = \F_{113}.
\end{equation}
\begin{equation}
\QS_{11|3} = \overline{\G_{131}}=\overline{\G_{311}}\,, \qquad
\QS_{113} = \overline{\G_{113}}.
\end{equation}
\begin{equation}
\QS_{12|34} = \overline{\G_{1324}} = \overline{\G_{3124}} =
\overline{\G_{1342}} = \overline{\G_{3142}} = \overline{\G_{3412}}\,.
\end{equation}

\medskip
\begin{theorem}
\label{QS-quot}
The $\PS_\qrpark$ form a basis of a Hopf subalgebra of $\PQSym$,
denoted by $\SQSym$. Its dual $\SQSym^*$ is the quotient
$\PQSym^*/{\mathcal J}$ where $\mathcal J$ is the two-sided ideal generated
by
\begin{equation}
\{\G_\park - \G_{\park'} | \park\congru\park' \}\,.
\end{equation}

Moreover, one has $\G_\park\congru \G_{\park'}$ iff
$\park\congru \park'$, so that $\SQSym^* \simeq \PQSym^*/\congru$.
The dual basis of $(\PS_\qrpark)$ is then $(\QS_{\qrpark})$.

The dimension of the component of degree $n$ of $\SQSym$ and $\SQSym^*$ is the
little Schr\"oder number (or super-Catalan) $s_n$.
\end{theorem}

\Proof
Let us begin with the elements $\overline\G_\park$. Since the
hypoplactic equivalence is a congruence
\begin{equation}
u\congru u' \text{\ and\ } v\congru v' \Longrightarrow uv\congru u'v',
\end{equation}
so that these elements build up an algebra that we will denote by $\SQSym^*$.
Since the hypoplactic congruence is compatible with the restriction to
intervals, one easily checks that the coproduct of $\G_\park$ is compatible
with the hypoplactic congruence, so that $\SQSym^*$ is a Hopf algebra.

Recall that two words $u$ and $v$ are hypoplactically equivalent iff they have
the same evaluation and $\Std(u)$ and $\Std(v)$ are hypoplactically
equivalent. Since two words of the same evaluation have parkized
words of the same evaluation as well, the same result applies if one replaces
the standardization by the parkization: two words $u$ and $v$ of the same
evaluation are hypoplactically equivalent iff their parkized words are.
This proves that
\begin{equation}
\park \congru \park' \Longleftrightarrow
\overline\G_\park=\overline\G_{\park'}.
\end{equation}
So $\SQSym^*$ is isomorphic to $\PQSym^*/\congru$ as a Hopf algebra.

Since the dual basis of $\G_\park$ in $\PQSym$ is $\F_\park$, one can write
the duality bracket as
\begin{equation}
\langle \G_\park, \F_{\park'} \rangle = \delta_{\park,\park'},
\end{equation}
where $\delta$ is the Kronecker symbol.
Then the dual basis of $\QS_q = \overline\G_\park$ inherited from the dual
Hopf algebras $\PQSym,\,\PQSym^*$ is naturally
$\sum_{\park'\congru\park}\F_{\park'}$, that is, $\PS_\qrpark$.
It then comes without proof that the $\PS_\qrpark$ form a basis of a Hopf
subalgebra of $\PQSym$ we will denote by $\SQSym$, as it is the dual of
$\SQSym^*$.

The dimensions are given by the little Schr\"oder numbers since these numbers
count the hypoplactic classes of parking functions. 
\qed

\begin{theorem}
The product and coproduct rules for the $\QS_\qrpark$ and the $\PS_\qrpark$
are

\begin{equation}
\QS_{\qrpark'} \QS_{\qrpark''} = 
\sum_{\park \in \park'\convP\park''} \QS_{\overline\park}\,,
\end{equation}
where $\park'$ (resp. $\park''$) is in the hypoplactic class of $\qrpark'$
(resp. $\qrpark''$).
\begin{equation}
\Delta \QS_{\qrpark} = \sum_{u,v ; \qrpark= u|v[|u|] \text{\ or\ }
                                   \qrpark= u.v[|u|]}
                   {\G_{u} \otimes \G_{v}}\,.
\end{equation}
\begin{equation}
\PS_{\qrpark'} \PS_{\qrpark''} = \PS_{\qrpark'|\qrparkb''} +
\PS_{\qrpark'\qrparkb''}
\end{equation}
where $\qrparkb''=\qrpark''[|\qrpark'|]$.
\begin{equation}
\Delta\PS_{\qrpark} =
\sum_{\qrpark',\qrpark''}\PS_{\qrpark'}\otimes\PS_{\qrpark''}
\end{equation}
where the sum is taken over the hypoplactic classes $\qrpark'$ and $\qrpark''$
such that their canonical elements $c'$ and $c''$ can be obtained as parkized
words of the prefix and the suffix of an element of the hypoplactic class
$\qrpark$.
\end{theorem}

\Proof
The formulas for the product and coproduct of the $\QS$ come from the formulas
of the $\G$ in $\PQSym^*$. The formulas for the $\PS$ are then easily derived
from the previous ones by duality.
\qed

For example,
\begin{equation}
\QS_{1|2}\QS_{1} = \QS_{1|23} + \QS_{1|22} + \QS_{12|3} + \QS_{11|3} +
\QS_{11|2} + \QS_{1|2|3}.
\end{equation}
\begin{equation}
\begin{split}
\Delta \QS_{11|34|55} =& 1\otimes\QS_{11|34|55} + \QS_{11}\otimes\QS_{12|33} +
\QS_{11|3}\otimes\QS_{1|22} \\ &+ \QS_{11|34}\otimes\QS_{11} +
\QS_{11|34|55}\otimes1.
\end{split}
\end{equation}
\begin{equation}
\PS_{11|335|6}\PS_{112} = \PS_{11|335|6778} + \PS_{11|335|6|778}.
\end{equation}
\begin{equation}
\Delta\PS_{11|3} = 1\otimes\PS_{11|3} +
\PS_{1}\otimes(\PS_{1|2}+\PS_{11}) +
(\PS_{21}+\PS_{12})\otimes\PS_{1} + \PS_{11|3}\otimes1.
\end{equation}

\subsection{$\SQSym$ is not self-dual}

Some simple computations prove that $\SQSym$ and $\SQSym^*$ are not isomorphic
Hopf algebras since the primitive Lie algebra of $\SQSym^*$ is of dimension
$6$ in degree $3$, spanned by
\begin{equation}
\begin{split}
& \QS_{111};\ \QS_{112};\ \QS_{11|2};\ \QS_{122}-\QS_{1|22};\
\QS_{113}-\QS_{11|3}; \\
&\QS_{123}-\QS_{1|23}-\QS_{12|3}+\QS_{1|2|3},
\end{split}
\end{equation}
whereas it is of dimension $7$ in $\SQSym$, spanned by:
\begin{equation}
\begin{split}
&\PS_{123} - \PS_{112} - \PS_{122} + \PS_{111};\ \
\PS_{1|22} - \PS_{11|2} - \PS_{112} + \PS_{122};\ \
\\ &
\PS_{1|23} - \PS_{11|2} - \PS_{112} + \PS_{111};\ \
\PS_{12|3} - \PS_{11|2} - \PS_{112} + \PS_{111};\ \
\\ &
\PS_{113} - \PS_{112};\ \
\PS_{11|3} - \PS_{11|2};\ \
\PS_{1|2|3} -\PS_{11|2} + \PS_{1|22}.
\end{split}
\end{equation}

In particular, it is impossible to endow $\SQSym$ or $\SQSym^*$ with a
bidendriform bialgebra structure since both would then be self-dual.
We cannot use the machinery of Foissy to investigate the freeness of both
algebras and their primitive Lie algebras, but we can do it by hand.

\subsection{Algebraic structure of $\SQSym^*$ and $\SQSym$}

Since we know that the primitive Lie algebra of $\SQSym$ is of dimension seven
in degree $3$, $\SQSym^*$ cannot be free and, indeed, one finds the relation
\begin{equation}
\QS_1 (\QS_{11}+\QS_{12}) = (\QS_{11} + \QS_{12}) \QS_1.
\end{equation}

We now move to $\SQSym$.
Consider the set $\PQS$ of parking quasi-ribbons that cannot be obtained as a
nontrivial shifted concatenation of parking quasi-ribbons. They are
the parking quasi-ribbons having a bar whenever the underlying nondecreasing
parking function has a breakpoint. For example, here are the elements of
$\PQS_n$ for $n\leq4$.
\begin{equation}
\begin{split}
&
 \{1\} ;\ \ \{11,\, 1\sepb2\} ;\ \
 \{111,\, 112,\, 11\sepb2,\, 11\sepb3,\, 1\sepb22,\, 1\sepb2\sepb3\} ;\\
& \{1111,\,1112,\,111\sepb2,\, 1113,\, 111\sepb3,\, 111\sepb4,\, 1122,\,
11\sepb22,\\
&\ 1123,\, 11\sepb23,\, 112\sepb3,\, 11\sepb2\sepb3,\, 112\sepb4,\,
 11\sepb2\sepb4,\, 11\sepb33, \\
&\ 11\sepb3\sepb4,\, 1\sepb222,\, 1\sepb223,\, 1\sepb22\sepb3,\,
 1\sepb22\sepb4,\, 1\sepb2\sepb33,\, 1\sepb2\sepb3\sepb4\}.
\end{split}
\end{equation}

Since the elements of $\PQS$ are those that never occur in a nontrivial
shifted concatenation of elements of $\PQR$, any element $\qrpark$ of $\PQR$
decomposes uniquely as a shifted product
$\qrpark_1\sconc \qrpark_2 \sconc\cdots\sconc \qrpark_k$ where all the
$\qrpark_k$ are in $\PQS$.
Define then
\begin{equation}
\PS^{\qrpark} = \PS_{\qrpark_1}\PS_{\qrpark_2}\ldots \PS_{\qrpark_k}.
\end{equation}

\begin{proposition}
The $\PS^\qrpark$ form a multiplicative basis of $\SQSym$.
In particular, $\SQSym$ is free as an algebra.
\end{proposition}

\Proof
The $\PS^\qrpark$ generate the same algebra as the $\PS_\qrpark$ since they
are triangular over the $\PS_\qrpark$: each term $\PS^{\qrpark}$ begins with
$\PS_{\qrpark}$ followed with elements of $\PQR$ that are shifted
concatenations of strictly lower elements of $\PQS$.
\qed

Since $\SQSym$ is free, one can compute the generating series of its
generating set. Recall that the generating series of $s_n$ is
$S(t):=\frac{1+t-\sqrt{1-6t+t^2}}{4t}$, so that

\begin{equation}
U(t):=1-1/S(t) = \frac{1-t-\sqrt{1-6t+t^2}}{2},
\end{equation}
that is the generating series of large Schr\"oder numbers $s'_n$ (sequence
A006318~\cite{Slo}), obviously equal to $2s_{n}$ thanks to the previous
formula. So
\begin{proposition}
The sets $\PQS_n$ are enumerated by the large Schr\"oder numbers.
\end{proposition}

\Proof
Even if the algebraic construction has already proved this result, we provide
a bijective proof in order to enlighten the relation between the large and the
little Schr\"oder numbers from the point of view of parking functions.

We split $\PQS_n$ into two and provide a bijection between both sets and
$\PQR_{n-1}$, the set of parking quasi-ribbons of length $n-1$.

Let $\PQS'_n$ be the subset of $\PQS_n$ consisting of the elements whose
underlying parking function is prime. The bijection between $\PQS'_n$ and
$\PQR_{n-1}$ is trivial: it consists in adding or removing $1$ at the
beginning of the parking function.

Let $\PQS''_n$ be the complementary subset of $\PQS_n$.
The bijection is the following: start from an element of $\PQR_{n-1}$. If this
element belongs to $\PQS_{n-1}$, then add a bar and $n$ to its end. Otherwise,
let $i$ be the smallest integer greater than $1$ such that $i-1$ is a
breakpoint and that there is no bar before the first $i$.
Then insert a bar and an $i$ before the first $i$.
This element satisfies the requirements of $\PQS_n$ since it can have
breakpoints only to the left of $i$ and that, by hypothesis, all those
breakpoints followed by a bar. Moreover, this element has a breakpoint,
so belongs to $\PQS''_n$. 
For example, the image of $11\sepb2\sepb455\sepb669$ is
$11\sepb2\sepb4\sepb55\sepb669$, since there is a breakpoint at $4$ with no
bar before the first $5$.

The reverse bijection consists in considering the rightmost breakpoint $i$ of
the underlying parking function of an element of $\PQS''_n$ and remove $i+1$
with the bar before it. The result belongs to $\PQR_{n-1}$ since we
removed the letter just after the \emph{rightmost} breakpoint.

Finally, it is a bijection between $\PQS''_n$ and $\PQR_{n-1}$ since the
operations are inverse to each other and the image of each set is included in
the other.
\qed

The next proposition summarizes the structures of $\SQSym$ and $\SQSym^*$.
\begin{proposition}
The algebra $\SQSym$ is a Hopf algebra of dimension $s_n$ in degree $n$.
It is not self-dual since $\SQSym$ is free as an algebra whereas $\SQSym^*$ is
not. 
\end{proposition}

\subsection{$\SQSym^*$ as a combinatorial Hopf algebra}

The embedding of Formula~(\ref{phi-isom}) induces an embedding
\begin{equation}
\QSym \simeq \FQSym^*/({\mathcal J}\cap\FQSym^*) \rightarrow
\PQSym^*/{\mathcal J} = \SQSym^*\,.
\end{equation}
In particular, we see that $\SQSym^*$ contains a large commutative subalgebra.

\subsection{Primitive Lie algebras of $\SQSym$ and $\SQSym^*$}

Since $\SQSym^*$ contains $\QSym$ as a subalgebra, its primitive Lie algebra
cannot be free and one easily finds:
\begin{equation}
[ \QS_{1}, \QS_{12} - \QS_{1|2} + \QS_{11}] =0.
\end{equation}

The first dimensions for the primitive Lie algebra of $\SQSym$ are
$1,\ 2,\ 7,\ 25,\ 102,\ $ with no relations between those elements in those
degrees so that one can conjecture that it is free as a Lie algebra.

\subsection{Schr\"oder ribbons}

In the algebra $\NCSF$, the product of non-commutative complete fonctions
split into sums of ribbon Schur functions, using a simple order on
compositions.
To get an analogous construction in our case, we define a partial order on
segmented non-decreasing parking functions.

Let $\pi$ be a segmented non-decreasing parking function and $\ev(\pi)$ be its
segmented evaluation vector, that is its evaluation vector with separators
between the $i$-th and $i+1$-th element if $i$ and $i+1$ are separated by a
bar in $\pi$. The successors of $\pi$ are the segmented non-decreasing
parking functions whose evaluations are given by the following algorithm:
given two non-zero elements of $\ev(\pi)$ not separated by a bar with only
zeroes between them, replace the left one by the sum of both and the right one
by 0.

For example, the successors of $11|3346$ are $11|3336$, and $11|3344$.

By transitive closure, the successor map gives rise to a partial order
$\succeq$ on segmented non-decreasing parking functions.

Now, define the Schr\"oder ribbons by
\begin{equation}
\label{SchRub}
\PS_\qrpark =: \sum_{\qrpark'\succeq\qrpark} {\RS_{\qrpark'}}\,.
\end{equation}
or, by M\"obius inversion on the boolean lattice,
\begin{equation}
\label{invSchRub}
\RS_\qrpark := \sum_{\qrpark'\succeq\qrpark} {(-1)^{f(\qrpark,\qrpark')}
\PS_{\qrpark'}}\,,
\end{equation}
where $f(u,v)$ is the difference between the numbers of different letters
in $u$ and in $v$.

For example,
\begin{equation}
\PS_{11|34} = \RS_{11|34} + \RS_{11|33}.
\end{equation}
\begin{equation}
\PS_{11|3346} = \RS_{11|3346} + \RS_{11|3336} + \RS_{11|3344} + \RS_{11|3333}.
\end{equation}
\begin{equation}
\RS_{11|3346} = \PS_{11|3346} - \PS_{11|3336} - \PS_{11|3344} + \PS_{11|3333}.
\end{equation}

\begin{proposition}
The product of two ribbons is given by
\begin{equation}
\label{prodRS}
\RS_{\qrpark'} \RS_{\qrpark''} = \RS_{\qrpark'|\qrparkb''} +
\RS_{\qrpark'\qrparkb''} + \RS_{\qrpark'\triangleright\qrparkb''}
\end{equation}
where $\qrparkb''=\qrpark''[|\qrpark'|]$ and
$\qrpark'\triangleright\qrparkb''$ is the successor of
$\qrpark'\qrparkb''$ obtained by decreasing the smallest letters of
$\qrparkb''$ down to the value of the greatest letters of $\qrpark'$.
\end{proposition}

\Proof
Let $p$ be the length of $\qrpark'$.
Let us expand $\RS_{\qrpark'}\RS_{\qrpark''}$ on the $\PS$ basis.
On gets an alternating sum of $\PS$ indexed by the successors of
$\qrpark'\qrparkb''$ having different $p$-th and $p+1$-st letters or indexed
by the successors of $\qrpark'|\qrparkb'$. This second set obviously sum up
to $\RS_{\qrpark'|\qrparkb''}$.
This first set is part of all successors of $\qrpark'\qrparkb''$, the missing
set being all successors of $\qrpark'\triangleright\qrparkb''$.

The sign of an element depending only on its number of different letters,
the result follows.
\qed

For example,
\begin{equation}
\RS_{1}\RS_{1|2} = \RS_{1|2|3} + \RS_{12|3} + \RS_{11|3}.
\end{equation}
\begin{equation}
\RS_{11|3}\RS_{113} = \RS_{11|3|446} + \RS_{11|3446} + \RS_{11|3336}.
\end{equation}

\subsection{Dendriform structures on $\SQSym$}

Let us now consider the other structures that can be put on $\SQSym$ and
$\SQSym^*$. First note that the product rules of $\PQSym^*$ as a tridendriform
algebra are compatible with the hypoplactic congruence, so that $\SQSym^*$ is
a tridendriform algebra. But it is not free since
\begin{equation}
\QS_{11|2} = \overline\G_{212} = \overline\G_{221},
\end{equation}
that can be rewritten as
\begin{equation}
(\QS_1\droittrid \QS_1)\miltrid\QS_1 = \QS_1\miltrid (\QS_1\droittrid \QS_1),
\end{equation}
a relation that is not a consequence of the tridendriform relations.

We already mentioned that $\SQSym^*$ cannot have a bidendriform bialgebra
structure since it would imply that $\SQSym^*$ is self-dual. On our
realization of the bidendriform bialgebra $\PQSym$, the explanation
comes from the fact that the hypoplactic congruence is not compatible with the
codendriform definitions since, for example,
\begin{equation}
\Delta_\gaudend \G_{221} = \G_1\otimes\G_{11} \text{\ whereas\ }
\Delta_\gaudend \G_{212} = 0.
\end{equation}

\newpage
\section{The Catalan Quasi-Symmetric Hopf algebra $\CQSym$}
\label{cqsym}

\subsection{The Hopf algebra $\CQSym$}

\subsubsection{Non-decreasing parking functions and non-crossing partitions}

As already mentioned, non-decreasing parking functions form a Catalan set.
There are dozens of possibilities to identify them to other combinatorial
objects. However, parking functions are known to be related to non-crossing
partitions (see~\cite{Bi1,St,Stan2}), and there is a simple bijection between
non-decreasing parking functions and non-crossing partitions. Starting with a
non-crossing partition, \emph{e.g.},
\begin{equation}
\pi = 13 | 2 | 45\,,
\end{equation}
one replaces all the letters of each block by its minimum, and reorders them
as a non-decreasing word
\begin{equation}
13 | 2 | 45 \to 11244,
\end{equation}
which is a parking function. In the sequel, we identify non-decreasing parking
functions and non-crossing partitions via this bijection.

\subsubsection{The Catalan Hopf algebra $\CQSym$}

For a general $\park\in\PF_n$, let $\NC(\park)$ be the non-crossing partition
corresponding to $\park^\uparrow$ by the inverse bijection, \emph{e.g.},
$\NC(42141)=\pi$ as above.

Then define $\PCat^\pi$ as the sum of all permutations of the non-decreasing
word corresponding to the given non-crossing partition:
\begin{equation}
\PCat^\pi := \sum_{\park ; \NC(\park)=\pi} {\F_\park}.
\end{equation}

\begin{theorem}
The $\PCat^\pi$, when $\pi$ runs over non-crossing partitions span a
cocommu\-ta\-tive Hopf subalgebra of $\PQSym$ with product and coproduct
given by
\begin{equation}
\label{prodP}
\PCat^{\pi'} \PCat^{\pi''} = \PCat^{\pi' \sconc \pi''}\,.
\end{equation}
\begin{equation}
\label{coprodP}
\Delta\PCat^{\pi} = \sum_{u,v ; (u.v)^\uparrow=\pi}
{\PCat^{\Park(u)} \otimes \PCat^{\Park(v)}}\,,
\end{equation}
where $u$ and $v$ run over the set of non-decreasing words.

Moreover, as an algebra, it is isomorphic to the algebra of the free semigroup
of non-crossing partitions under the operation of concatenation of diagrams.
\end{theorem}

\begin{proof}
Equation~(\ref{prodP}) follows from Equation~(\ref{prodF}): indeed,
any permutation of $\pi' \sconc \pi''$ is uniquely obtained as the shifted
shuffle of a permutation of $\pi'$ with a permutation of $\pi''$. The converse
is obvious.

Equation~(\ref{coprodP}) comes from Equation~(\ref{coprodF}):
consider the relation $P(p,q)$ on words which consists of pairs $(w,w')$ of
words $w$ and $w'$ of length $p+q$ such that the sorted word of the prefix of
length $p$ (resp. suffix of length $q$) of $w$ and $w'$ are equal.
By definition of $\PCat^\pi$, it is a sum of such classes, so that
$\Delta \PCat^\pi$ decomposes as a sum of tensor products of the form
$\PCat^{\pi'}\otimes\PCat^{\pi''}$. The sum on the right hand-side of
Equation~(\ref{coprodP}) is exactly over representatives of the equivalence
classes, hence the result.
Formula~(\ref{coprodP}) proves that the coalgebra $\CQSym$ is cocommutative.

Moreover, since $\CQSym$ is a subalgebra and a sub-coalgebra of $\PQSym$, the
product and the coproduct of $\CQSym$ are compatible, so that $\CQSym$ is
endowed with a graded bialgebra structure, and therefore, with a Hopf algebra
structure.
\end{proof}

This algebra will be called the \emph{Catalan subalgebra} of $\PQSym$ and
denoted by $\CQSym$.

For example, one has
\begin{equation}
\PCat^{11} \PCat^{1233} = \PCat^{113455} \qquad
\PCat^{1124} \PCat^{1223} = \PCat^{11245667}
\end{equation}
\begin{equation}
\begin{split}
\Delta\PCat^{1124} &= 1\otimes\PCat^{1124} +
        \PCat^1\otimes \left(\PCat^{112}+\PCat^{113}+\PCat^{123}\right) +
        \PCat^{11}\otimes\PCat^{12} \\
       & + \PCat^{12}\otimes\left(\PCat^{11}+2\PCat^{12}\right)
         + \left(\PCat^{112}+\PCat^{113}+\PCat^{123}\right)\otimes\PCat^1 +
        \PCat^{1124}\otimes1\,.
\end{split}
\end{equation}

Since the non-decreasing parking functions that never occur in a nontrivial
shifted concatenation of such elements are the connected non-decreasing
parking functions, any $\pi$ decomposes uniquely as a shifted product
$\pi_1\sconc \pi_2 \sconc\cdots\sconc \pi_k$ where all the
$\pi_k$ are connected.

\begin{proposition}
The $\PCat$ form a multiplicative basis of $\CQSym$.
In particular, $\CQSym$ is free as an algebra.
\end{proposition}

Here are the connected non-decreasing parking functions up to
length $4$.
\begin{equation}
\{1\}, \{11\}, \{111, 112\}, \{1111, 1112, 1113, 1122, 1123\}.
\end{equation}

Since $\CQSym$ is free, one can compute the generating series of its
generating set. Recall that the generating series of $C_n$ is
$C(t):=\frac{1-\sqrt{1-4t}}{2t}$, so that
\begin{equation}
CN(t):=1-1/C(t) = t \frac{1-\sqrt{1-6t+t^2}}{2},
\end{equation}
that is the generating series of shifted Catalan numbers $C_{n-1}$.
Indeed, the connected non-decreasing parking functions are obtained by
concatenating a $1$ to the left of all non-decreasing parking functions.

\subsubsection{Algebraic structure of $\CQSym$}

Following Reutenauer~\cite{Re} p.~58, denote by $\pi_1$ the Eulerian
idempotent, that is, the endomorphism of $\CQSym$ defined by
$\pi_1=\log^*(Id)$ where $\log^*$ means that the logarithm is taken in the
convolution algebra of graded endomorphisms $\End^{\rm gr}(\CQSym)$.
It is obvious, thanks to the definition of $\PCat^\pi$ that
\begin{equation}
\pi_1(\PCat^\pi) = \PCat^\pi + \cdots,
\end{equation}
where the dots stand for terms $\PCat^\gamma$ where $\gamma$ is not connected.
So the family $\pi_1(\PCat^\alpha)$ where $\alpha$ runs over all connected
non-decreasing parking functions is a free set of primitive generators of
$\CQSym$. In particular, they generate a free Lie algebra (see,
\emph{e.g.}~\cite{HNT3} for more details) whose Hilber series is given by
\begin{equation}
t + t^2 + 3\, t^3 + 8\,t^4 + 25\,t^5 + 75\,t^6 + 245\, t^7 + 800\,t^8+
O(t^9).
\end{equation}
The sequence is referenced in Sloane's database as A022553~\cite{Slo}. It
counts Lyndon words $l$ of even length $2n$ with an equal number of $a$ and
$b$
So the free Lie algebra of primitive elements of $\CQSym$ is isomorphic to the
Lie subalgebra $\L$ of the free Lie algebra
$\Lie(a,b)$ consisting of the elements with an equal number of
$a$ and $b$.
One can then prove that the standard bracketings of
the Lyndon words $l$ with the same  number of $a$ and $b$
such that $l=l'\cdot b$ with $l'$ also being a Lyndon word generate a free Lie
algebra. Since those particular Lyndon words are enumerated by the shifted
Catalan numbers, one can conclude that they generate $\L$.

\subsection{The dual Hopf algebra $\CQSym^*$}

Let us denote by $\MM_\pi$ the dual basis of $\PCat^\pi$ in the
\emph{commutative} algebra $\CQSym^*$. Since $\CQSym$ is the subalgebra of
$\PQSym$ obtained by summing all permutations of nondecreasing parking
functions, $\CQSym^*$ is the quotient of $\PQSym^*$ by the relations
$\G_\a\equiv\G_\b$ if $\a^\uparrow=\b^\uparrow$.

It is then immediate (see Equation~(\ref{prodG})) that the multiplication in
this basis is given by
\begin{equation}
\label{prodMM}
\MM_{\pi'} \MM_{\pi''} = \sum_{
\a\in\pi'\convol\pi''}{\MM_{\a^\uparrow}}\,.
\end{equation}

For example,
\begin{equation}
\MM_{1} \MM_{12} = \MM_{112} + \MM_{113} + \MM_{122} + 3 \MM_{123}.
\end{equation}
\begin{equation}
\begin{split}
\MM_{12} \MM_{11} &= \MM_{1112} + \MM_{1113} + \MM_{1114} + \MM_{1123} +
\MM_{1124}\\
&+ \MM_{1134} + \MM_{1222} + \MM_{1223} + \MM_{1224} + \MM_{1233}\,.
\end{split}
\end{equation}

\begin{theorem}
$\CQSym^*$ can be embedded in the polynomial algebra $\CC[x_1,x_2,\ldots]$ by
\begin{equation}
\MM_{\pi} = \sum_{\Park(w)=\pi} {\underline{w}}\,,
\end{equation}
where $\underline{w}$ is the commutative image of $w$
(\emph{i.e.}, $a_i\mapsto x_i$).
\end{theorem}

\begin{proof}
The result follows  from (\ref{realG}) and from the fact that
quotienting  $\PQSym^*$ by the relations $\G_\park\equiv\G_\parkb$ if
$\park^\uparrow=\parkb^\uparrow$  amounts to take the commutative image
of words, transforming these into monomials.
\end{proof}

For example,
\begin{equation}
\MM_{111} = \sum_{i} {x_i^3}\,.
\end{equation}
\begin{equation}
\MM_{112} = \sum_{i} {x_i^2x_{i+1}}\,.
\end{equation}
\begin{equation}
\MM_{113} = \sum_{i,j ; j\geq i+2} {x_i^2 x_j}\,.
\end{equation}
\begin{equation}
\MM_{122} = \sum_{i,j ; i<j} {x_i x_j^2}\,.
\end{equation}
\begin{equation}
\MM_{123} = \sum_{i,j,k; i<j<k} {x_ix_jx_k}\,.
\end{equation}

The \emph{packed evaluation vector} $\tasse(w)$ of $w$ is obtained
from $\ev(w)$ by removing all its zeroes.
For example, if $w=3117291781329$, $\ev(w)=(4,2,2,0,0,0,2,1,2)$ and
$\tasse(w)=(4,2,2,2,1,2)$.

We can now see that $\CQSym^*$  contains $\QSym$ as a subalgebra.
The embedding of $\QSym$ into $\CQSym^*$ is given by
\begin{equation}
\label{gammaEmb}
\gamma(M_I) := \sum_{\tasse(\pi)=I} {\MM_{\pi}}.
\end{equation}

For example,
\begin{equation}
M_{3} = \MM_{111},\quad
M_{21}= \MM_{112}+\MM_{113},\quad
M_{12}= \MM_{122},\quad
M_{111}=\MM_{123}.
\end{equation}

\subsection{Catalan ribbons}

As already done for the Schr\"oder algebras, we define a partial order on
non-decreasing parking functions.

Let $\pi$ be a non-decreasing parking function and $\ev(\pi)$ be its
evaluation vector. The successors of $\pi$ are the non-decreasing parking
functions whose evaluations are given by the following algorithm: given two
non-zero elements of $\ev(\pi)$ with only zeroes between them, replace the
left one by the sum of both and the right one by 0.

For example, the successors of $113346$ are $111146$, $113336$, and $113344$.

By transitive closure, the successor map gives rise to a partial order on
non-decreasing parking functions. We will write $\pi'\succeq\pi$ if $\pi'$ is
obtained from $\pi$ by successive applications of successor maps.

Now, define the Catalan ribbon functions by
\begin{equation}
\label{catalRub}
\PCat^\pi =: \sum_{\pi'\succeq\pi} {\RCat_{\pi'}}\,.
\end{equation}
This last equation completely defines the $\RCat_{\pi}$.

For example,
\begin{equation}
\begin{split}
\PCat^{113346} = &\ \ \
\RCat_{113346} + \RCat_{113344} + \RCat_{113336} + \RCat_{113333}\\
& + \RCat_{111146} + \RCat_{111144} + \RCat_{111116} + \RCat_{111111}
\end{split}
\end{equation}
and
\begin{equation}
\begin{split}
\RCat_{113346} = &\ \ \
\PCat^{113346} - \PCat^{113344} - \PCat^{113336} + \PCat^{113333}\\
& - \PCat^{111146} + \PCat^{111144} + \PCat^{111116} - \PCat^{111111}.
\end{split}
\end{equation}
Note that, by M\"obius inversion of the boolean lattice, the coefficient of
$\PCat^{\pi'}$ in $\RCat^{\pi}$ is $-1$ to the number of different letters in
$\pi$ minus the number of different letters in $\pi'$.

This definition is compatible with the definition of commutative ribbon Schur
functions since if one considers the morphism $\phi$ as
\begin{equation}
\begin{split}
\phi:
\CQSym &\mapsto Sym \\
\PCat^\pi &\to S^{c(\pi)}
\end{split}
\end{equation}
then the image $\phi(\RCat_{\pi})$ is equal to $R_{c(\pi)}$.

\begin{proposition}
The product of two $\RCat$ functions is
\begin{equation}
\label{prodRCat}
\RCat_{\pi'} \RCat_{\pi''} = \RCat_{\pi'\sconc\pi''} +
\RCat_{\pi'\triangleright\pi''}\,,
\end{equation}
where $\triangleright$ is the successor of $\pi'\sconc\pi''$ obtained by
decreasing the smallest letters of $\pi''$ down to the greatest letters of
$\pi'$.
\end{proposition}

\begin{proof}
Let $p$ be the length of $\pi'$.
Let us expand $\RCat_{\pi'\sconc\pi''}$ on the $\PCat$ basis.
On gets the alternating sum of $\PCat$ indexed by successors of
$\pi'\sconc\pi''$. Those successors split into two disjoint subsets: the
successors having the $p$-th and $p+1$-th letters equal and the others.
The first set corresponds to the successors of $\pi'\triangleright\pi''$
whereas the second set corresponds to the $w'\sconc w''$ where
$w'\succeq \pi'$ and $w''\succeq \pi''$.

The sign of an element depending only on its number of different letters,
the alternating sum of the first set amounts to
$-\RCat_{\pi'\triangleright\pi''}$ whereas the sum of the second set amounts
to $\RCat_{\pi'} \RCat_{\pi''}$.
\end{proof}

For example,
\begin{equation}
\RCat_{11224} \RCat_{113} = \RCat_{11224668} + \RCat_{11224448}\,,\quad
\RCat_{113} \RCat_{11224} = \RCat_{11344557} + \RCat_{11333557}\,.
\end{equation}

\subsection{Internal product}

Define the parkized word of a bimonomial as the non-decreasing parking
function obtained by parkizing its lexicographically sorted biword.
Recall that bimonomials can be encoded as matrices, the entry $A_{ij}$
being the number of bi-letters $(ij)$ in the biword, so that it makes sense to
speak of the parkized word of a matrix.

\begin{theorem}[~\cite{NT2}]
\label{cqsym-stab}
The homogeneous components $\CQSym_n$ of the Catalan algebra are
stable under the internal product $*$.
More precisely, one has
\begin{equation}
\label{prodintP}
\P^{\pi'}*\P^{\pi''} = \sum_{\pi} \P^\pi
\end{equation}
where $\pi$ runs over the parkized words of all non-negative integer matrices
with row sum $\ev(\pi')$ and column sum $\ev(\pi'')$.
\end{theorem}

\begin{example}
\begin{equation}
\P^{1123} * \P^{1111} = \P^{1134}; \qquad
\P^{1111} * \P^{1123} = \P^{1123}.
\end{equation}
\begin{equation}
\P^{1123} * \P^{1112} = 2\P^{1134} + \P^{1234}; \qquad
\P^{1122} * \P^{1224} = \P^{1134} + \P^{1233} + 2\P^{1234}.
\end{equation}
\begin{equation}
\P^{1123} * \P^{1224} = 2\P^{1134} + 5\P^{1234}.
\end{equation}
\end{example}

The matrices appearing in the last product are
\begin{equation}
\begin{split}
\begin{pmatrix}
1&1&.&.\\ .&1&.&. \\ .&.&.&1
\end{pmatrix}
\begin{pmatrix}
1&1&.&.\\ .&.&.&1 \\ .&1&.&.
\end{pmatrix}
\begin{pmatrix}
1&.&.&1\\ .&1&.&. \\ .&1&.&.
\end{pmatrix}
\begin{pmatrix}
.&2&.&.\\ 1&.&.&. \\ .&.&.&1
\end{pmatrix}
\\[6pt]
\begin{pmatrix}
.&2&.&.\\ .&.&.&1 \\ 1&.&.&.
\end{pmatrix}
\begin{pmatrix}
.&1&.&1\\ 1&.&.&. \\ .&1&.&.
\end{pmatrix}
\begin{pmatrix}
.&1&.&1\\ .&1&.&. \\ 1&.&.&.
\end{pmatrix}
\end{split}
\end{equation}
the fourth and the fifth matrices having $1134$ as parkized word whereas the
other ones yield $1234$.

It is interesting to observe that these algebras are non-unital.
Indeed, it follows from Formula~(\ref{prodintP}) that

\begin{corollary}
The element $\J_n=\P^{(1^n)}$ is a left unit for $*$, but not
a right unit.
\qed
\end{corollary}

The description of $\P^{\pi'}*\P^{\pi''}$ in terms of integer matrices being
essentially identical to that of $S^I*S^J$ in $\NCSF$, the same argument as
in~\cite{NCSF1}, proof of Proposition 5.2, shows that the splitting formula
remains valid in $\CQSym_n$:

\begin{proposition}
\label{splitting}
Let $\mu_r$ denote the $r$-fold product map from $\CQSym^{\otimes r}$
to $\CQSym$, $\Delta^r$ the $r$-fold coproduct with values in
$\CQSym^{\otimes r}$, and $*_r$ the internal product of the $r$-fold tensor
product of algebras $\CQSym^{\otimes r}$. Then, for
$f_1,\ldots,f_r,g\in\CQSym$,
\begin{equation}
(f_1\cdots f_r)*g=\mu_r[(f_1\otimes\cdots\otimes f_r)*_r \Delta^r(g)]\,.
\hskip6.3cm
\qed
\end{equation}
\end{proposition}

This is indeed the same formula as with the internal product of $\Sym$,
actually, an extension of it, since we have

\begin{corollary}
The Hopf subalgebra of $\CQSym$ generated by the elements $\J_n$, which is
isomorphic to $\Sym$ by $j:S_n\mapsto \J_n$, is stable under $*$, and thus
also $*$-isomophic to $\Sym$.
Moreover, the map $f\mapsto f*\J_n$ is a projector onto $\Sym_n$, which is
therefore a left $*$-ideal of $\CQSym_n$.
\qed
\end{corollary}

More precisely, if $i<j<\ldots<r$ are the letters occuring in $\pi$, so that
as a word $\pi=i^{m_i}j^{m_j}\cdots r^{m_r}$, then
\begin{equation}
\P^\pi*\J_n = \J_{m_i}\J_{m_j}\cdots \J_{m_r}\,.
\end{equation}

It follows from Theorem~\ref{cqsym-stab} that the $\RCat_{\pi}$ are the
pre-images of the ordinary ribbons under the projection $f\mapsto f*\J_n$:
\begin{corollary}
Let $I$ be the composition obtained by discarding the zeros of the evaluation
of an non-decreasing parking function $\pi$. Then
\begin{equation}
\RCat_\pi * \J_n = j(R_I).
\end{equation}
More precisely, if $I=(i_1,\ldots,i_p)$, this last element is equal to
$\RCat_{1^{i_1}\sconc1^{i_2}\sconc\cdots\sconc1^{i_p}}$, that is, the Catalan
ribbon indexed by the only non-decreasing word of evaluation $d(\pi)$.
\qed
\end{corollary}

The internal product of $\CQSym$ is dual to the coproduct $\delta f=f(XY)$ on
the commutative algebra $CQSym$, quotient of $\PQSym^*$.
For example, we have
\begin{equation}
\begin{split}
\MA_{113}(XY) =& (\MA_{112}(X)+\MA_{113}(X))
                (\MA_{111}(Y)+\MA_{112}(Y)+\MA_{113}(Y)+\MA_{122}(Y)) \\
&+ \MA_{111}(X) \MA_{113}(Y).
\end{split}
\end{equation}
\begin{equation}
\MA_{112}(XY) =  \MA_{111}(X) \MA_{112}(Y).
\end{equation}

\subsubsection{Cauchy Kernel}

Define the \emph{Cauchy kernel} by
\begin{equation}
K(X;A)=\sum_{\a\in\PF}\G_\a(X)\F_\a(A)=\sum_\pi \MM_\pi(X)\P^\pi(A)\,.
\end{equation}

\begin{proposition}
The kernel $K$ has the reproducing property
\begin{equation}
K(X;A)*K(Y;A)=K(XY;A)\,.
\end{equation}
\end{proposition}

\Proof
\begin{equation}
\begin{split}
\<K(X)*K(Y),\MM_\pi\>=
 & \sum_{\pi',\pi''}\MM_{\pi'}(X)\MM_{\pi''}(Y)
     \<\P^{\pi'}*\P^{\pi''},\MM_\pi\>\\
=& \sum_{\pi',\pi''}\MM_{\pi'}(X)\MM_{\pi''}(Y)
     \<\P^{\pi'}\otimes\P^{\pi''} ,\Delta\MM_\pi\>\\
=& \Delta\MM_\pi(X,Y)=\MM_\pi(XY).
\end{split}
\end{equation}
\qed

\subsection{Compositions, Lagrange inversion, and $H_n(0)$}
\label{compos-sec}

\subsubsection{}

Recall that nondecreasing parking functions (or non-crossing partitions) can
be classified according to the factorization
$\pi = \pi_1\sconc \cdots\sconc \pi_r$ into irreducible nondecreasing parking
functions (or non-crossing partitions). Let $i_k:=|\pi_k|$ and
$c(\pi):=(i_1,\cdots,i_k)$, regarded as a composition of $n$.

We set
\begin{equation}
\V^I := \sum_{c(\pi)=I} \PCat^{\pi}
\end{equation}
considered as an element of $\PQSym$.
If one defines $\V_n = \V^{(n)}$, we have
\begin{equation}
\V_n = \sum_{\park\in\PPF_{n}} {\F_\park}
\end{equation}
and
\begin{equation}
\V^I = \V_{i_1} \cdots \V_{i_r} = \sum_{\park\in\PPF_I} {\F_\park}\,.
\end{equation}
This can be reformulated as
\begin{equation}
\sum_{\a\in\PF}\F_\a = \left(1-\sum_{\b\in\PPF}\F_\b\right)^{-1}\,,
\end{equation}
which is the lift to $\PQSym$ of the well known identity
\begin{equation}
\sum_{n\geq 0}(n+1)^{n-1}\frac{t^n}{n!}=
\left(1-\sum_{n\geq 1}(n-1)^{n-1}\frac{t^n}{n!}\right)^{-1}\,.
\end{equation}
Indeed, the map $\F_\a\mapsto \frac1{n!}$ ($n=|\a|$) is a character
of $\PQSym$.

At this point, it is useful to observe that if $C(w)$ denotes the descent
composition of a word $w$, the map
\begin{equation}
\eta: \F_\park \mapsto F_{C(\park)}\,,
\end{equation}
which is a Hopf algebra morphism $\PQSym\to\QSym$, maps $\V^I$ to the
Frobenius characteristic of the underlying permutation representation of
$\SG_n$ on $\PPF_I$.
\begin{equation}
\eta(\V^I) = \sum_{\park\in\PPF_I}\F_{C(\park)} = \text{ch}(\PPF_I)\,,
\end{equation}
Indeed, if $V\subseteq A^n$ is any set of words invariant under the
right action of $\SG_n$, the characteristic of the underlying permutation
representation is always equal to $\sum_{w\in V}F_{C(w)}$. This is because
$V$ splits as a disjoint union
\begin{equation}
V=\bigsqcup_\nu A_\nu\,,\ \text{where}\ A_\nu=\{w\in A^n|\ev(w)=\nu\}\,.
\end{equation}
The characteristic of $A_\nu$ is clearly $h_\nu$, and it is well
known that
\begin{equation}
h_\nu=\sum_{w\in A_\nu}F_{C(w)}\,.
\end{equation}
Actually, each $\C A_\nu$ is also a projective $H_n(0)$-module
with noncommutative characteristic ${\bf ch}(\C A_\nu)=S^I$,
where $I=\tasse(\nu)$.

\subsubsection{}

As a consequence, the number of parking functions of type $I$ with descent
composition $J$ is equal to the scalar product of symmetric functions
\begin{equation}
\< r_J, f^I\>
\end{equation}
where $f^I=f_{i_1}\cdots f_{i_r}=\text{ch}(\PPF_I)$ and $r_J$ is the ribbon
Schur function. This extends Prop.~3.2.(a) of~\cite{St}.
Remark that in particular, by inversion of 
\begin{equation}
\F_{\PF_{n}} := \sum_{\park\in\PF_{n}} {\F_\park} =
  \sum_{I\vDash n} {\V^I}\,,
\end{equation}
 one obtains
\begin{equation}
\F_{\PPF_n} = \sum_{I\vDash n} {(-1)^{n-l(I)}\F_{\PF_I}}\,,
\end{equation}
where
\begin{equation}
\PF_I := \PF_{i_1} \ssh \PF_{i_2} \ssh \cdots \ssh \PF_{i_r}\,.
\end{equation}

These identities are easily visualized on the encoding of parking functions
with skew Young diagrams as in~\cite{PP} or in~\cite{HHLRU}.

\subsubsection{}

The transpose $\gamma^*$ of the map $\gamma$ defined in
Equation~(\ref{gammaEmb}), is the map
\begin{equation}
\begin{split}
{\bf ch } :& \CQSym \to \NCSF \\
           & \PCat^\pi \mapsto  S^{\tasse(\pi)}\,.
\end{split}
\end{equation}
which sends $\PCat^\pi$ to the characteristic non-commutative symmetric
function of the natural projective $H_n(0)$-module with basis
$\{\park\in\PF_n|\NC(\park)=\pi\}$.

\subsubsection{}

One can show that
\begin{equation}\label{defg}
g := \sum_{n\geq0} g_n := \sum_{n\geq0} {\bf ch}(\F_{\PF_n}) =
\sum_{I} {\bf ch}(\V^I).
\end{equation}
is the series obtained by applying the non-commutative Lagrange inversion
formula of~\cite{Ges,PPR} to the generating series of complete functions,
\emph{i.e.}, $g$ is the unique solution of the equation
\begin{equation}\label{lagr}
g = 1 + S_1g + S_2g^2 + \cdots = \sum_{n\geq0}{S_n g^n}\,.
\end{equation}
Indeed, let $g$ be defined by (\ref{defg}) and set
\begin{equation}
f := \sum_{n\geq 1} f_n= \sum_{n\geq 1}{\bf ch}(\F_{\PPF_n})=\sum_{n\geq 1}{\bf
ch}(\V_n)\,.
\end{equation}
Recall that the prime nondecreasing parking functions of length $n$ are
obtained by concatenating a 1 to the left of a nondecreasing parking function
of length $n-1$. This gives a recurrence for $f_n$ and $g_n$.
From each nondecreasing parking function of packed evaluation
$J\vDash n-1$, contributing a term $S^J$ to $g_{n-1}$, 
we get a prime nondecreasing parking function of
packed evaluation $I=(j_1+1,j_2,\ldots,j_r)\vDash n$, contributing
a term $S^I :=\Omega S^J$ to $f_n$, where $\Omega$ is the linear
operator incrementing the first part in the basis $S^J$ of $\Sym$.
Hence, $f_n=\Omega g_{n-1}$, and we have a system of two equations
\begin{equation}\label{eqsfg}
\left\{
\begin{array}{lcl}
f &=& \Omega g\\
g &=& (1-f)^{-1}\,,
\end{array}
\right.
\end{equation}
which, with the initial condition $g_0=1$, admits a unique solution:
$f_1=\Omega g_0=S^1$, $g_1=f_1$, $f_2=\Omega g_1=S^2$, 
$g_2=f_2+f^{11}=S^2+S^{11}$, $f_3=\Omega g_2=S^3+S^{21}$,
$g_3=f_3+f^{21}+f^{12}+f^{111}=S^3+2S^{21}+S^{12}+S^{111}$, and so on.

But the unique solution of (\ref{lagr}) satisfies
\begin{equation}
\Omega g= S_1+S_2g+S_3g^2+\cdots
\end{equation}
and also
\begin{equation}
1=g^{-1}+ S_1+S_2g+S_3g^2+\cdots = g^{-1}+\Omega g
\end{equation}
so that if we set $f=\Omega g$, we solve (\ref{eqsfg}) as well.

Remark that the commutative images of these equations give  the
$\SG_n$-characteristics of $\PF_n$ and $\PPF_n$, and that we have derived them
from first principles, using only the multiplication rule of $\PQSym$ and the
notion of a prime parking function.

\subsubsection{The Hilbert series of $\SQSym$ revisited}

It is also possible to obtain it by a character calculation,
derived from the above considerations.
If we decompose the noncommutative characteristic of the $H_n(0)$-module
$\C\PF_n$ into ribbons
\begin{equation}
{\bf ch}(\C\PF_n)=\sum_{I\vDash n} m_I R_I\,,
\end{equation}
the number of hypoplactic classes of parking functions of length $n$
is
\begin{equation}
\sum_{I\vDash n} m_I\,.
\end{equation} 

Indeed, as already mentionned, if $V\subset A^n$ is any set of words which is
a disjoint union of evaluation classes $A_\nu$, $\CC V$ is a projective
$H_n(0)$-module since it is the direct sum $\oplus\CC A_\nu$, where
${\bf ch}(\CC A_\nu)=S^I$, with $I=\tasse(\nu)$.

Now, each $A_\nu$ is itself a disjoint union of hypoplactic classes
\begin{equation}
A_{\nu,I} = \{ w\in A_\nu | C(\Std(w)^{-1}) = I\},
\end{equation}
and each such class is the support of an indecomposable projective module
\begin{equation}
{\bf ch}(A_{\nu,I}) = R_I.
\end{equation}
By duality between the bases $F_I$ and $R_I$,
\begin{equation}
\sum{m_I} = \left\langle \sum F_I, {\bf ch}(\CC\PF_n) \right\rangle
\end{equation}
and taking into account the identity~\cite{Ung}
\begin{equation}
\sum_I F_I = \frac12 \left[1+\prod_{i\ge 1}\frac{1+x_i}{1-x_i}\right]\,
= 1 + \frac12 \sum_{n\geq1} \sum_{k=0}^n e_k h_{n-k},
\end{equation}
we obtain
\begin{equation}
\begin{split}
\dim(\SQSym_n) &= \left\< \sum_{I\vDash n} F_I, {\bf ch}(\F_{\PF_n}) \right\>
= \left\< \frac{1}{2}\sum_{k=0}^n{e_kh_{n-k}}, \frac{1}{n+1}h_n((n+1)X)
\right\> \\
&= \frac{1}{2n+2} \sum_{k=0}^n{\binom{n+1}{k}\binom{2n-k}{n-k}} = s_n\,.
\end{split}
\end{equation}

\newpage
\section{A Hopf algebra of segmented compositions}
\label{cc}

\subsection{Segmented compositions}

Define a \emph{segmented composition} as a finite sequence of positive
integers, separated by vertical bars or commas, \emph{e.g.},
$(2,1\sep2\sep1,2)$.

The number of segmented compositions having the same underlying composition is
obviously $2^{l-1}$ where $l$ is the length of the composition, so that the
total number of segmented compositions of sum $n$ is $3^{n-1}$ since
$(1+2)^{n-1} = 3^{n-1}$.

\subsection{A Hopf subalgebra of $\SQSym^*$}

\subsubsection{Hypoplactic packed words}

Let $A=\{a_1<a_2<\ldots\}$ be an infinite totally ordered alphabet. The
\emph{packed word} $u=\pack(w)$ associated with a word $w\in A^*$ is
obtained by the following process. If $b_1<b_2<\ldots <b_r$ are the letters
occuring in $w$, $u$ is the image of $w$ by the homomorphism
$b_i\mapsto a_i$.

A word $u$ is said to be \emph{packed} if $\pack(u)=u$. We denote by $\PW$ the
set of packed words.

Let us consider packed quasi-ribbons, that is, quasi-ribbons that are packed
words.
For example, the word $11324355$ is a packed quasi-ribbon word since it is the
reading of the following quasi-ribbon
\begin{equation}
\label{ex-packqr}
\PetitTableau
\Tableau{1&1&2\\\ &\ &3&3\\\ &\ &\ &4&5&5\\}
\end{equation}
These objects are in bijection with segmented compositions. Indeed, start from
a packed quasi-ribbon $\qrpark$ and write the evaluation vector $\Ig$ of
$\qrpark$, putting a separator between $\Ig_i$ and $\Ig_{i+1}$ iff $i$ and
$i+1$ are not is the same row of $\qrpark$. For example, the segmented
composition corresponding to the quasi-ribbon of Equation~(\ref{ex-packqr}) is
$21|2|12$. This element will be denoted by $\psev(\qrpark)$.
The reverse bijection consists in writing the unique nondecreasing word
of evaluation $\Ig$ and put the letters $i+1$ on the row next to the row
of letters of $i$ iff $i$ and $i+1$ are separated by $\sep$ in $\Ig$.

\begin{example}
For $n=2$, we have $3$ packed quasi-ribbons
\begin{equation}
11
\quad
12 
\quad
1|2 
\end{equation}

For $n=3$, we have $9$ packed quasi-ribbons
\begin{equation}
111\ \ 112\ \ 11\sep2 \ \ \ 122\ \ 1\sep22\ \
123\ \ 1\sep23\ \ 12\sep3\ \ 1\sep2\sep3.
\end{equation}
respectively encoded as the $9$ segmented compositions:
\begin{equation}
3\ \ 21\ \ 2|1\ \ 12\ \ 1|2\ \ 111\ \ 1|11\ \ 11|1\ \ 1|1|1.
\end{equation}
\end{example}

In the sequel, we will identify packed quasi-ribbons and their encodings as
segmented compositions.

\subsubsection{A Hopf subalgebra of $\SQSym^*$}

Let us denote by ${\sf P}(w)$ the hypoplactic $P$-symbol of a word $w$ (its
quasi-ribbon). The $P$-symbols of packed words are therefore packed
quasi-ribbons.

For each packed quasi-ribbon $\Ig$, define
\begin{equation}
\PSC_{\Ig} := \PS_{\detassmax(\Ig)}
\quad\text{and}\quad
\QSC_{\Ig} := \sum_{\psev(\qrpark)= \Ig} \QS_\qrpark \in \SQSym^*.
\end{equation}
where $\detassmax(\qrpark)$ is the maximal parking quasi-ribbon for the
lexicographic order of evaluation $\Ig$.

For example,
\begin{equation}
\PSC_{11|2} = \PS_{12|33}, \qquad
\PSC_{112} = \PS_{1233}.
\end{equation}
\begin{equation}
\QSC_{12|1} = \QS_{122|3} + \QS_{122|4}, \qquad
\QSC_{121} = \QS_{1223} + \QS_{1224}.
\end{equation}
\begin{equation}
\QSC_{12|21} = \QS_{122|334} + \QS_{122|335} + \QS_{122|336} + \QS_{122|445}
+ \QS_{122|446}.
\end{equation}

\begin{theorem}
The $\PSC_\Ig$ span a Hopf subalgebra $\SCQSym$ of $\SQSym$. This subalgebra
is also the quotient of $\SCQSym$ by the relations
$\PS_\qrpark=\PS_{\qrpark'}$ if $w$ and $w'$ have same packed word.

The $\QSC_\Ig$ span a Hopf subalgebra $\SCQSym^*$ of $\SQSym^*$ and one has
\begin{equation}
\dim\ \SCQSym_n=3^{n-1},
\end{equation}
for $n\geq1$.
\end{theorem}

\Proof
The product and coproduct rules of the $\PS$ of $\SQSym$ imply that
the $\PSC$ of $\SQSym^*$ span a Hopf subalgebra of $\SQSym$. It is also
obvious that both operations are compatible with the relations
$\PS_\qrpark=\PS_{\qrpark'}$ if $\qrpark$ and $\qrpark'$ have same packed
word, so that $\SCQSym$ is a Hopf quotient of $\SQSym$. Then, by the usual
standard argument, we have that the $\QSC$ are a basis of $\SCQSym^*$.

The dimension of $\SCQSym_n$ is given by the number of segmented compositions,
that is $3^{n-1}$.
\qed

To describe the product and coproduct rules of both bases, we need a new
operation on segmented compositions.

Recall that the product of the monomial basis $M_I$ on $\QSym$ is defined with
the augmented shuffle of two compositions recursively defined as
\begin{equation}
(I_1,I') \assh (J_1,J') =
I_1 (I'\assh (J_1,J')) + J_1 ((I_1,I')\assh J') + I_1+J_1 (I'\assh J'),
\end{equation}
with the extra condition $I'\assh\epsilon=\epsilon\assh I'$ where $\epsilon$
is the empty word.

This construction is generalized to $\SCQSym$ as follows: the \emph{augmented
shuffle} $\Ig'\assh \Ig''$ of two segmented compositions is obtained from the
usual augmented shuffle $I'\assh I''$ of their underlying compositions by
inserting bars between two blocks $I_k$ and $I_{k+1}$ of a composition $I$ iff
\begin{itemize}
\item $I_k$ and $I_{k+1}$ both contains elements coming from $\Ig'$ and those
elements were separated by a bar,
\item $I_k$ and $I_{k+1}$ both contains elements coming from $\Ig''$ and those
elements were separated by a bar,
\item $I_k$ contains an element coming from $\Ig''$ and $I_{k+1}$ contains an
element coming from $\Ig'$.
\end{itemize}

For example,
\begin{equation}
1 \assh 2|1 = 12|1 + 3|1+ 2|11 + 2|2 + 2|1|1, \quad
1 \assh 21  = 121  + 31 + 2|11 + 2|2 + 21|1.
\end{equation}

\begin{theorem}
The product and coproduct rules for the $\PSC_\Ig$ and the $\QSC_\Ig$
are
\begin{equation}
\PSC_{\Ig'} \PSC_{\Ig''} = \PSC_{\Ig'|\Ig''} + \PSC_{\Ig'\Ig''}.
\end{equation}
\begin{equation}
\Delta\PSC_{\Ig} =
\sum_{\Ig \in \Ig'\assh\Ig''}\PS_{\Ig'}\otimes\PS_{\Ig''}.
\end{equation}
\begin{equation}
\label{QSC-prod}
\QSC_{\Ig'} \QSC_{\Ig''} = \sum_{\Ig\in \Ig'\assh\Ig''} \QSC_\Ig.
\end{equation}
\begin{equation}
\Delta \QSC_\Ig = \sum_{\Ig=\Ig'\cdot \Ig'' \text{\ or\ } \Ig=\Ig'|\Ig''}
  \QSC_{\Ig'}\otimes\QSC_{\Ig''}.
\end{equation}
\end{theorem}

\Proof
The product of two $\PSC$ of $\SCQSym$ directly comes from the product of two
$\PS$ of $\SQSym$. The coproduct of a $\QSC$ then comes by duality.
The shifted shuffle of two compositions then obviously give all the
possible evaluations of the convolution of two parking functions of the given
evaluations. Finally, the rules to place the bars correspond to the different
cases where there is an $i$ to the right of an $i+1$ in one of the resulting
parking functions.
\qed

For example,
\begin{equation}
\PSC_{12|1} \PSC_{2|11} = \PSC_{12|12|11} + \PSC_{12|1|2|11}.
\end{equation}
\begin{equation}
\begin{split}
\Delta\PSC_{12|1} =&\ \ \ \
1\otimes\PSC_{12|1} + \PSC_{1}\otimes (\PSC_{12}+\PSC_{2|1})
+ \PSC_{11}\otimes (\PSC_{1|1}+\PSC_{2}) \\ &+ \PSC_{1|1}\otimes\PSC_{2}
+ (\PSC_{111}+\PSC_{11|1}) \otimes \PSC_{1}
+ \PSC_{12|1}\otimes1.
\end{split}
\end{equation}
\begin{equation}
\QSC_{1}\QSC_{2|1} = \QSC_{3|1} + \QSC_{12|1} + \QSC_{2|2} + \QSC_{2|1|1}
+ \QSC_{2|11}.
\end{equation}
\begin{equation}
\QSC_{1}\QSC_{11|1} = \QSC_{111|1} + \QSC_{21|1} + \QSC_{1|11|1} +
\QSC_{1|2|1} + \QSC_{11|11} + \QSC_{11|2} + \QSC_{11|1|1}.
\end{equation}
\begin{equation}
\begin{split}
\QSC_{1|1}\QSC_{1|1} =&\ \ \ \
2 \QSC_{1|11|1} + \QSC_{1|1|11} + \QSC_{11|11} + \QSC_{11|1|1} +
\QSC_{1|1|1|1} \\ &+ \QSC_{2|11} + \QSC_{2|1|1} + \QSC_{2|2} + 2\QSC_{1|2|1} +
\QSC_{1|1|2} + \QSC_{11|2}.
\end{split}
\end{equation}
\begin{equation}
\Delta \QSC_{12|1} = 1\otimes\QSC_{12|1} + \QSC_{1}\otimes\QSC_{2|1} +
\QSC_{12}\otimes\QSC_{1} + \QSC_{12|1}\otimes1.
\end{equation}

\subsection{Algebraic structure of $\SCQSym$ and $\SCQSym^*$}

The algebra $\SCQSym^*$ is not free for exactly the same reason $\SQSym^*$ is
not: one has the relation
\begin{equation}
\QSC_1 (\QSC_{2}+\QSC_{11}) = (\QSC_{2} + \QSC_{11}) \QSC_1.
\end{equation}

Let us now move to $\SCQSym$.

Since $\SCQSym$ is the subalgebra of $\SQSym$ spanned by the parking
quasi-ribbons that are maximally unpacked, and since $\SQSym$ is free,
$\SCQSym$ is automatically free and generated by the maximal elements of
$\PQS$. For example, the generators of $\SCQSym$ for $n\leq4$ are
\begin{equation}
\begin{split}
&
 \{1\} ;\ \ \{11,\, 1\sepb2\} ;\ \
 \{111,\, 11\sepb3,\, 1\sepb22,\, 1\sepb2\sepb3\} ;\\
& \{1111,\, 111\sepb4,\, 11\sepb33,\, 11\sepb3\sepb4, 1\sepb222,\, 
 1\sepb22\sepb4,\, 1\sepb2\sepb33,\, 1\sepb2\sepb3\sepb4\},
\end{split}
\end{equation}
that can be rewritten on segmented compositions as
\begin{equation}
\begin{split}
&
 \{1\} ;\ \ \{2,\, 1\sepb1\} ;\ \
 \{3,\, 2\sepb1,\, 1\sepb2,\, 1\sepb1\sepb1\} ;\\
& \{4,\, 3\sepb1,\, 2\sepb2,\, 2\sepb1\sepb1, 1\sepb3,\, 
 1\sepb2\sepb1,\, 1\sepb1\sepb2,\, 1\sepb1\sepb1\sepb1\},
\end{split}
\end{equation}
By the same argument on generating series as in $\SQSym$, one finds that there
are $2^{n-1}$ generators of $\SCQSym$ of degree $n$. And indeed, these
generators are in natural bijection with compositions of $n$ since they
have separators between all elements.

The next proposition summarizes the structures $\SCQSym$ and $\SCQSym^*$.
\begin{proposition}
The algebra $\SCQSym$ is a Hopf algebra of dimension $3^{n-1}$.
It is not self-dual since $\SCQSym$ is free as an algebra whereas $\SCQSym^*$
is not.
Moreover, $\SCQSym$ is free over a graded alphabet labelled by all
compositions.
\end{proposition}

\subsection{Primitive Lie algebras of $\SCQSym^*$}

Since $\SCQSym^*$ contains $\QSym$ as a subalgebra, its primitive Lie algebra
cannot be free and one easily finds:
\begin{equation}
[ \QSC_{1}, \QSC_{12} - \QSC_{1|2} + \QSC_{11}] =0.
\end{equation}

\subsection{A quasi-ribbon basis of $\SCQSym^*$}

The elements $\QSC_\Ig$ are segmented analogs of the basis
$(M_I)$ of $\QSym$. So we can define analogs of the $(F_I)$ of $\QSym$ in the
same way as we did in $\SQSym^*$.

Recall that the \emph{refinement order} denoted by $\raff$ on compositions
is such that $I=(i_1,\ldots,i_k)\raff J=(j_1,\ldots,j_l)$ iff
$\{i_1,i_1+i_2,\ldots,i_1+\cdots+i_k\}$ contains
$\{j_1,j_1+j_2,\ldots,j_1+\cdots+j_l\}$.
In this case, we say that $I$ is finer than $J$.
For example, $(2,1,2,3,1,2)\raff (3,2,6)$.

Let $\Ig=(I_1\sep \cdots\sep I_r)$ and let

\begin{equation}
\FSC_\Ig := \sum_{\Ig'} \QSC_{\Ig'},
\end{equation}
where the sum is taken over sequences of compositions $(I'_1,\ldots,I'_r)$
where $I'_k$ is finer than $I_k$.
For example, one has
\begin{equation}
\FSC_{2|2} = \QSC_{11|11} + \QSC_{2|1} + \QSC_{11|2} + \QSC_{2|2}\,.
\end{equation}

By a triangularity argument, we have

\begin{theorem}
The $\FSC_\Ig$ form a basis of $\SCQSym^*$.
\end{theorem}

The basis $\FSC_\Ig$ satisfies a product formula similar to the $F_I$ of
$\QSym$ (whence the choice of notation).
To state it, we need an analogue of the shifted shuffle.

A \emph{segmented permutation} is a permutation with separators.
The \emph{descent composition} $C(\alpha)$ of a segmented permutation $\alpha$
is a segmented composition: it is the sequence of descent compositions of the
blocks of $\alpha$ separated by bars.

For example, $\alpha = 248|517|3$ is a segmented permutation whose descent
composition is $(3|12|1)$.

The \emph{shifted shuffle} $\alpha\ssh \beta$ of two segmented permutations is
obtained from the usual shifted shuffle $\sigma\ssh\tau$ of the underlying
permutations $\sigma$ and $\tau$ by inserting bars
\begin{itemize}
\item after each descent which was originally followed by a bar in $\alpha$
or in the shift of $\beta$,
\item after each descent created by the shuffling process.
\end{itemize}

For example,
\begin{equation}
2|1 \ssh 21 = 2|143 + 24|13 + 243|1 + 4|2|13 + 4|23|1 + 43|2|1.
\end{equation}

\begin{theorem}
Let $\Ig'$ and $\Ig''$ be two segmented compositions and let $\alpha$ and
$\beta$ be any two segmented permutations whose descent compositions are
respectively $\Ig'$ and $\Ig''$.
Then
\begin{equation}
\FSC_{\Ig'} \FSC_{\Ig''} = \sum_{\Ig} \FSC_\Ig,
\end{equation}
where the sum runs over the descent compositions of the segmented permutations
$\gamma$ occuring in $\alpha\ssh\beta$.
\end{theorem}

\begin{proof}
The product of the $\QSC$ can be easily rewritten in terms of of segmented
permutations as follows: associate with two segmented compositions $\Ig'$ and
$\Ig''$ two segmented permutations $\alpha$ and $\beta$ such that
$C(\alpha)=\Ig'$ and $C(\beta)=\Ig''$. Consider the elements in the shifted
shuffle $\alpha\ssh\beta$ such that two elements in increasing order of
$\alpha$ (resp. $\beta$) not separated by a bar have no $\beta$ (resp.
$\alpha$) between them.
From all those elements, build the set of all segmented permutations with at
least those bars and at most new bars between the elements of $\alpha$ and the
elements of $\beta$. For all those segmented permutations, compute first their
descent compositions and \emph{then} remove the bars added lately.
The set of the descent compositions obtained by this process correspond to the
product $\QSC_\alpha\QSC_\beta$.

Express both $\FSC$ in the $\QSC$ basis and group the terms in their product
where the letters of $\alpha$ have been inserted at the same place. By
construction, the lexicographically minimum element $s$ in each group with the
smallest number of bars belongs to $\alpha\ssh\beta$.
Now, given the product rule of the $\QSC$, we have all elements obtained from
$s$ by adding any number of bars, thus $\FSC_s$.
\end{proof}

For example,
\begin{equation}
\FSC_{1}\FSC_{11|1} = \FSC_{21|1} +\FSC_{1|2|1} +\FSC_{11|2} + \FSC_{11|1|1}.
\end{equation}
\begin{equation}
\FSC_{1} \FSC_{2|1} = \FSC_{3|1} + \FSC_{1|2|1} + \FSC_{2|2} + \FSC_{2|1|1}.
\end{equation}

\begin{theorem}
Let $\Ig$ be a segmented composition.
Then
\begin{equation}
\Delta \FSC_\Ig = \sum_{\Ig=\Ig'\cdot \Ig'' \text{\ or\ } \Ig=\Ig'|\Ig''
                         \text{\ or\ } \Ig=\Ig'\triangleright\Ig''}
  \FSC_{\Ig'}\otimes\FSC_{\Ig''},
\end{equation}
where $(I_1,\ldots,I_k)\triangleright (I'_1,\ldots,I'_l)$ denotes the
segmented composition $(I_1,\ldots,I_{k-1},I_k\triangleright'
I'_1,I'_2,\ldots,I'_l)$ where
$(i_1,\ldots,i_k)\triangleright' (j_1,\ldots,j_l) =
(i_1,\ldots,i_{k-1},i_k+j_1,j_2,\ldots,j_l)$.
\end{theorem}

\Proof
This result will follow by duality from the considerations in the forthcoming
section.
\qed

\subsection{A ribbon basis of $\SCQSym$}

Let $(\RSC_\Ig)$ be the dual basis of $(\FSC_\Ig)$.
Then

\begin{proposition}
The $(\RSC_\Ig)$ are a basis of $\SCQSym$ related to the $\PSC$ by
\begin{equation}
\PSC_\Ig =: \sum_{\Ig'} \RSC_{\Ig'},
\end{equation}
where the sum is taken over sequences of segmented compositions
$(I'_1|\ldots|I'_r)$ where $I_k$ is finer than $I'_k$.
\end{proposition}

Since $\SCQSym$ is a subalgebra of $\SQSym$ such that the image of $\PSC_\Ig$
is $\PS_\qrpark$, given both orders on parking quasi-ribbons and on segmented
compositions, the image of $\RSC_{\Ig}$ is $\RS_\qrpark$.
This remark immediately proves the product rule of the $\RSC$, whereas its
coproduct rule comes from the duality between the $(\RSC)$ and the $(\FSC)$.

\begin{theorem}
The product and coproduct rules of $\RSC$ are
\begin{equation}
\label{prodRIg}
\RSC_{\Ig'} \cdot  \RSC_{\Ig''} = \RSC_{\Ig'\cdot \Ig''} + \RSC_{\Ig'|\Ig''}
 + \RSC_{\Ig'\triangleright\Ig''}.
\end{equation}
\begin{equation}
\Delta \RSC_{\Ig} = \sum_{\Ig\in\Ig'\ssh\Ig''} \RSC_{\Ig'}\otimes\RSC_{\Ig''}.
\end{equation}
\end{theorem}

From Equation~(\ref{prodRIg}), we see that $\SCQSym$ is 
the free cubical trialgebra on one generator, see~\cite{LRtri}.

\newpage
\section{Appendix}
\label{app}

\subsection{Relations with free probability theory}

The free cumulants $R_n$ of a probability measure $\mu$ on $\RR$
are defined (see {\it e.g.,}~\cite{Spei}) by means of the generating series
of its moments $M_n$
\begin{equation}
G_\mu(z) :=\int_{\RR}\frac{\mu(dx)}{z-x}=z^{-1}+\sum_{n\ge 1}M_nz^{-n-1}
\end{equation}
as the coefficients of its compositional inverse 
\begin{equation}
K_\mu(z) :=G_\mu(z)^{\<-1\>}=z^{-1}+\sum_{n\ge 1}R_n z^{n-1}\,.
\end{equation}
It is in general instructive to interpret the coefficients of
a formal power series as the specializations of the elements of some 
generating family of the algebra of symmetric functions.
In this context, it is the interpretation
\begin{equation}
M_n=\phi(h_n)=h_n(A)
\end{equation}
which is relevant.
Indeed, the process of functional inversion (Lagrange
inversion) admits a simple expression within this formalism
(see~\cite{Mcd}, ex.~24 p.~35). If the symmetric functions
$h_n^*$ are defined by the equations
\begin{equation}
u=tH(t) \ \Longleftrightarrow \ t=uH^*(u)
\end{equation}
where $H(t):=\sum_{n\ge 0}h_nt^n$, $H^*(u):=\sum_{n\ge 0}h_n^*u^n$,
then, using the $\lambda$-ring notation,
\begin{equation}
h_n^*(X)=\frac{1}{n+1}(-1)^{n}e_n((n+1)X) :=\frac{1}{n+1} [t^n] E(-t)^{n+1}
\end{equation}
where $E(t)$ is defined by $E(t)H(t)=1$.
This defines an involution $f\mapsto f^*$ of the ring of symmetric
functions.

Now, if one sets $M_n=h_n(A)$ as above, then
\begin{equation}
G_\mu(z)=z^{-1}H(z^{-1})=u 
\ \Longleftrightarrow \ 
z= K_\mu(u)=\frac{1}{u}E^*(-u) =
u^{-1}+\sum_{n\ge 1}(-1)^ne_n^* u^{n-1}\,.
\end{equation}
Hence,
\begin{equation}  
R_n=(-1)^ne_n^*(A)\,.
\end{equation}

\subsection{An exercise on permutation representations}

It follows immediately from the explicit formula (see~\cite{Mcd} p. 35)
\begin{equation}
-e_n^*=\frac{1}{n-1}\sum_{\lambda\vdash n}
\binom{n-1}{l(\lambda)} \binom{l(\lambda)}{m_1, m_2,\ldots,m_n} e_\lambda
\end{equation}
(where $\lambda=1^{m_1}2^{m_2}\cdots n^{m_n}$) that $-e_n^*$
is Schur positive. Clearly, $-e_n^*$ is the Frobenius characteristic
of a permutation representation $\Pi_n$, twisted by the sign character.
Let us set
\begin{equation}
(-1)^{(n-1)}R_n = -e_n^*=:\omega(f_n)
\end{equation}
so that
\begin{equation}
\label{fn-def}
f_n := \sum_{\lambda\vdash n}
\frac{1}{n-1}
\binom{n-1}{l(\lambda)} \binom{l(\lambda)}{m_1, m_2,\ldots,m_n} h_\lambda
\end{equation}
and $f_n$ is the character of $\Pi_n$.

The problem of constructing such a representation had been
raised by Kerov in 1995. 
We shall see that $\Pi_n$ corresponds to prime parking functions.
We note that our construction of $\Pi_n$ is merely a variation about
previously known results (see in particular~\cite{Len,PP}). However, since
this is this precise version of the question that led us to the Hopf algebra
of parking functions and some of its properties, we decided to include its
discussion.

\subsection{Solution of the exercise}

\begin{proposition}
The Frobenius characteristic of the permutation representation of $\PPF_n$ is
$f_n$.
\end{proposition}

\Proof
We first show that the number of nondecreasing prime parking functions whose
reordered evaluation is a given partition $\lambda$ is equal to

\begin{equation}
\frac{1}{n-1} \binom{n-1}{l(\lambda)} \binom{l(\lambda)}{m_1,m_2,\ldots,m_n}
\end{equation}
where $\lambda=1^{m_1}2^{m_2}\cdots n^{m_n}$ (see Equation~(\ref{fn-def})).
Indeed, this number corresponds to the number of ways of putting the
$\lambda_i$ over $n-1$ places in a circle. For such a placement $P$,
number in all possible clockwise ways the places of the circle
and consider the $n-1$ nondecreasing words $i^{c_i}$ where $c_i$ is the
content of place number $i$.
Then, thanks to~\cite{DM}, there is exactly one of those words that is a prime
parking function.

For example, with $\lambda=(3,2,2)$, there are $10$ possible circles. Consider
the circle where the $3$ is followed by one empty place, a $2$, two
empty places, and the last $2$. Then the six nondecreasing words are:
\begin{equation}
1113366, 2255666, 1144555, 3344466, 2233355, 1122244,
\end{equation}
Now, since all permutations of a nondecreasing prime parking function are
parking functions, the Frobenius characteristic of the permutation
representation of this set of words is $h_\lambda$.
It then easily comes that
\begin{equation}
\ch(\PPF_n) = f_{n}\,,
\end{equation}
so that $\Pi_n$ can be identified with $\PPF_n$, as claimed before.
\qed

\newpage

\end{document}